\documentclass[titlepage]{amsart}
\usepackage{amsaddr}
 \usepackage{color, soul}  
\usepackage{ulem}                                                  
\usepackage{amsfonts}
\usepackage{amsmath}
\usepackage{setspace}
\usepackage{graphicx}
\usepackage{psfrag}
\usepackage{graphics}
\usepackage{amssymb}
\usepackage{bm}
\usepackage{multicol}
\usepackage{amsmath, amssymb,mathabx,mathrsfs}
\usepackage[toc,page]{appendix}
\usepackage{comment}
\usepackage{epsfig}
\usepackage{graphicx, subfigure}
\newtheorem{thm}{Theorem}[section]
\newtheorem{cor}[thm]{Corollary}
\newtheorem{lem}[thm]{Lemma}

\theoremstyle{definition}

\theoremstyle{remark}

\numberwithin{equation}{section}
\theoremstyle{property}
\newtheorem{prpty}[thm]{Property}
\theoremstyle{assumption}
\newtheorem{assum}[thm]{Assumption}
\theoremstyle{summary}
\newtheorem{sumry}[thm]{Summary}

\begin{document}
\title[ ]{Cantor Set Structure of the Weak Stability Boundary for Infinitely Many Cycles in the Restricted Three-Body Problem}

\author{EDWARD BELBRUNO}
\address{Yeshiva University, Department of Mathematical Sciences, New York, NY 10016, USA }
\email{edward.belbruno@yu.edu}



\begin{abstract}
The geometry of the weak stability boundary region for the planar restricted three-body problem about the secondary mass point has been an open problem.  Previous studies have conjectured that it may have a fractal structure.  In this paper, this region is studied for infinitely many cycles about the secondary mass point, instead of a finite number studied  previously.  It is shown that in this case the boundary consists of a family of infinitely many Cantor sets and is thus fractal in nature. It is also shown that on two-dimensional surfaces of section,  it is the boundary of a region only having bounded cycling motion for infinitely many cycles, while the complement of this region generally has unbounded motion.  It is shown that that this shares many properties of a Mandelbrot set.    Its relationship to the non-existence of KAM tori is described, among many other properties.  Applications are discussed. 
\end{abstract}

\keywords{celestial mechanics, chaos, stability, restricted three-body problem, KAM,  Cantor set}

\maketitle

\section{Introduction}  \label{sec:Intro} 

\medskip

The planar circular restricted three-body problem is considered. This model goes back to the work of H. Poincar\'e \cite{Poincare:1899}.  In it, a particle $P$ of negligible mass moves on the same plane as two mass-points, $P_1, P_2$, as they make mutual circular orbits about their common center of mass, with constant angular frequency.  $P_1$ is assumed to be much more massive than $P_2$. For example, $P_1$ is the Sun and $P_2$ is the Earth.   

The motion of $P$ about $P_1$ is well understood by the Kolmogorov-Arnold-Moser (KAM) Theorem when $P$ starts with elliptic Keplerian initial conditions about $P_1$.  
When the ratio of the angular frequencies of the circular motion of $P_1,P_2$ to that of the initial elliptic orbit of $P$ are sufficiently irrational, satisfying diophantine conditions, then  
the motion of $P$ about $P_1$ remains close to elliptical for all time, and is quasi-periodic. It lies on two-dimensional invariant tori on a three-dimensional fixed energy surface in the four-dimensional phase space of position and velocity. The tori separate the energy surface into disjoint connected components. There are gaps in the tori corresponding to approximate rational frequency ratios where the motion of $P$ is chaotic \cite{SiegelMoser:1971, Moser:1962, Arnold:1989, Celletti:2007}.

The motion of $P$ about $P_2$ is a much different situation. 
When $P$ starts with Keplerian elliptic initial conditions about $P_2$, the resulting dynamics is not generally understood. If $P$ starts sufficiently near to $P_2$, where the Jacobi energy is sufficiently large, then KAM tori exist by a theorem of M. Kummer \cite{Kummer:1979}. This is accomplished by a special regularization. As the distance from $P_2$ increases, however,  this methodology is not possible and the existence of invariant tori is not proven to exist.  This is an open problem as to how far from $P_2$ invariant tori can exist, although there are interesting numerical analysis  (see \cite{Simo:2000}).

A method to understand the motion about $P_2$ is to use weak stability boundary methods.  The weak stability boundary, $W$,   defines a region about $P_2$ in phase space that measures the stability of cycling motion about $P_2$.  It was first defined in \cite{Belbruno:1987} in 1987 for the more general three-dimensional restricted three-body problem (see also \cite{Belbruno:2004, Belbruno:2007}). The motivation for this was to find special low energy trajectories from the Earth to the Moon that were dynamically captured at the Moon. This led to a new type of transfer to the Moon for spacecraft which use no fuel for lunar capture.\footnote{It was  demonstrated operationally in 1991 with Japan's Hiten spacecraft \cite{BelbrunoMiller:1993}. This transfer class has been used by other lunar missions: GRAIL, CAPSTONE, Flashlight (NASA),  Denuri (KARI), Hakuto-R (ispace)}  Related capture transfers, obtained by considering $W$, were used also for the  SMART-1 lunar mission of the European Space Agency (ESA) \cite{Schoenmaekers:2001} as well as ESA's BepiColombo mission to Mercury \cite{Jehn:2008} (see also, \cite{Circi:2001}, \cite{Romagnoli:2009}).

$W$ is determined numerically, by estimating a transition between {\it stable} and {\it unstable} motion after $n$ cycles of the trajectory for $P$ about $P_2$, $n = 1, 2, \ldots $,  starting from the set of Keplerian elliptical conditions at the periapsis positions  (see Section \ref{sec:W}).  This yields a distribution of points about $P_2$ depending on several parameters, including $n$, defining  the $n$th weak stability boundary, $W_n$.  This set is a boundary for stable cycling motion. The case of $n=1$ was defined in \cite{Belbruno:1987}.  A more general and rigorous analysis for any integer $n \geq 1$ was done by E. Garcia, G. Gomez in an important paper \cite{GarciaGomez:2007}.   This was further numerically and theoretically studied by M. Gidea, F. Topputo, E. Belbruno in \cite{BelbrunoGideaTopputo:2010}, with a general numerically based result showing the equivalence of $W_n$ to a network of intersecting invariant manifolds associated to Lyapunov orbits (see Section \ref{subsec:GlobalW}). $W$ does not exist on a fixed Jacobi energy surface, but rather for a range of energies.

The general structure of  $W$ has been an open problem. This paper determines the geometry of this set in the case of infinite cycling and shows it is an infinite family of Cantor sets. More precisely, this paper shows the equivalency of the limit of $W_n$ as $ n \rightarrow \infty$, labeled $W'$,  with the boundary of a set $M^*$, $W'= \partial M^*$, and that  $W'$  is a Cantor set of hyperbolic points on each two-dimensional Poincar\'e surface of section, $S_{\theta}$, for a Poincar\'e map $\Phi$, parameterized by the polar angle $\theta \in [0, 2\pi]$ with respect to $P_2$, for each fixed Jacobi energy, $C$, and mass ratio, within a suitable range. $W'$ is an infinite union of Cantor sets of hyperbolic points over all $\theta$, suitable Jacobi energies,  and mass ratios. {It is important to note that the key result of this paper is Lemma \ref{lem:CantorSet} on the Cantor structure of $W'$.  This lemma depends on three assumptions. These are Assumptions \ref{assum:1}, \ref{assum:2}, \ref{assum:3}. They are based on numerical results which support their validity, which is discussed.   This implies that this paper doesn't present a completely rigorous analytic proof of the Cantor structure for $W'$. In other parts of the paper rigorous proofs and non-rigorous proofs are distinguished when necessary for clarity. If assumptions are used, their validity is discussed if necessary and consequences if they are not satisfied.}

 The interior of $M^*$ corresponds to points yielding stable motion for infinitely many cycles. $W'$ represents points giving unstable motion for infinitely many cycles.  $W'$ is the weak stability boundary for infinitely many cycles of $P$ about $P_2$ (see Theorem \ref{thm:M*W'}).

 It is remarked that $M^*$ shares a number of properties with that of the classical Mandelbrot set even though they are defined much differently (see Table \ref{tab:Comparison}).   The main difference is that the boundary of  $M$ is a continuous fractal set, whereas the boundary of $M^*$ is totally discontinuous on each two-dimensional surface of section. 

It is also remarked that the existence of $W'$ being an infinite union of Cantor sets (see Section \ref{subsec:GlobalW}) 
clarifies a brief comment in \cite{GarciaGomez:2007}(page 3) alluding to the possibility the stable points themselves may be related to a Cantor set, which is not entirely correct since it turns out it is the unstable points (see Section \ref{subsec:Properties}).  

The relationship of $M^*$ to the existence and non-existence of KAM tori is discussed in Section \ref{sec:NoKAM}  (Theorem \ref{thm:NoKAM}).  It says that for points of $\mathcal{P}$  on the boundary of $M^*$ then these cannot lie on KAM tori, while points on the interior of $M^*$ may lie on KAM tori, but that is not known.

A general discussion of the results is given in Section \ref{sec:Discussion}. 

The methodology of this paper is to give a theoretical framework of the results, which is supported and motivated by numerous previous numerical studies.

\section{Restricted Three-Body Problem} 
\label{sec:RestrictedProblem}

\begin{figure}
\centering
          \includegraphics[width=0.65\textwidth, clip, keepaspectratio]{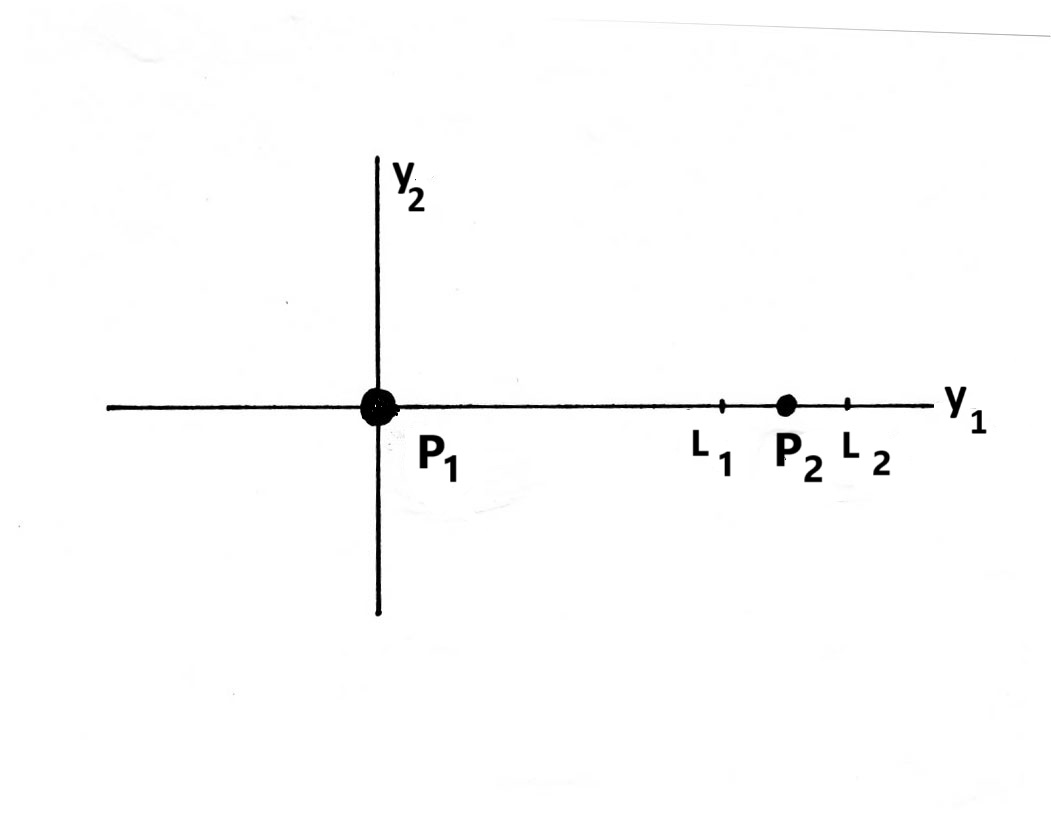}
	\caption{Rotating coordinate system $(y_1,y_2)$. }
	\label{fig:Coors}
\end{figure}

The two mass points $P_1, P_2$, defined in the Introduction,  are assumed to move in mutual circular orbits of constant frequency, $\omega = 1$, about their common center of mass on a plane in a rotating coordinate system $y_1, y_2$ of rotational frequency $\omega$, centered at $P_1$.   The mass of $P_1$ is $1-\mu$ and the mass of $P_2$ is $\mu$, $\mu$ is assumed to be small.  $P_1, P_2$ are fixed on the $y_1$-axis at  $(0,0), (1,0)$, respectively (see  Figure \ref{fig:Coors}).\footnote{The center of mass between $P_1, P_2$  is at $(\mu, 0)$.} The gravitational constant is normalized to $1$. The particle $P$ of zero mass moves on the $(y_1, y_2)$-plane. Its motion is given by ${\bf y}(t) = (y_1(t), y_2(t))$ as a function of time, $t$, defined as a solution to the system of differential equations,
 
   \begin{equation}
   \ddot{y}_1 - 2 \dot{y}_2  = \Omega_{y_1} ,  \hspace{.1in}
   \ddot{y}_2 +2 \dot{y}_1  = \Omega_{y_2} ,   \hspace{.1in},
       \label{eq:DEs}
   \end{equation}
  where 
  \begin{equation}
  \Omega = \frac{1}{2} [(y_1-\mu)^2  + y_2^2]  + \frac{1-\mu}{r_1} + \frac{\mu}{r_2} + \frac{1}{2} \mu(1-\mu),
  \label{eq:Omega}
  \end{equation}
  and $r_1^2 = y_1^2 +y_2^2 $, $ r_2 ^2 = (y_1 - 1)^2 + y_2^2 $, $\Omega_x\equiv \partial \Omega / \partial x$,  $^. \equiv$ d/dt.  
  \medskip

 The Jacobi integral of the system is
  \begin{equation}
  J = J({\bf y}, \dot{\bf{y}}  ) = 2 \Omega - |\dot{\bf{y}}|^2,  
  \label{eq:JacobiEnergy}
  \end{equation}
  $|\dot{\bf{y}}|^2 =  (\dot{y}_1^2 + \dot{y}_2^2 )$.
  \medskip
 Thus, along a solution $ {{\bm{\psi}}} (t) = ({\bf{y}}(t), \dot{\bf{y}}(t) )$,    
 \begin{equation}
 J({\bm{\psi}} (t))  = C,
 \label{eq:Jacobi}
 \end{equation} 
 for a constant $C$, the Jacobi constant.   Each value of $C$ defines a three-dimensional energy surface in the four-dimensional phase 
 space $({\bf{y}}, \dot{\bf{y}})$ that the solutions lie on,  the Jacobi surface, $\sigma = \sigma(C) = \{({\bf{y}}, \dot{\bf{y}}) \in \mathbb{R}^4| J=C\} $.  
   
  As $C$ varies, these surfaces have different geometric properties that constrain the motion of $P$ in ${\bf y}$-space. This is obtained by projecting the Jacobi surfaces onto ${\bf y}$-space, defining the Hill\rq{}s regions where physical motion can occur.  The Hill's regions are defined by ${H}(C) = \{ {\bf y} \in \mathbb{R}^2 | 2\Omega \geq C \}$.   $(\ref{eq:DEs})$ has five equilibrium points, $L_i, 1=1,2,3,4,5$, the Euler-Lagrange points, where $\dot{\bf{y}} = \ddot{\bf{y}} = \bf{0}$, where $C=C_i, C_4 = C_5 = 3 < C_3 < C_2 < C_1$.  $L_1, L_2$ are of interest to this paper and shown in Figure {\ref{fig:Coors}.

 The set $Z =\{ {\bf y}\in \mathbb{R}^2 | 2\Omega = C\}$ defines the points of zero velocity, which is the boundary curve of the Hill\rq{}s regions.   When $C > C_1$, trajectories cannot pass between $P_1$ and $P_2$ from the respective disconnected Hill's regions about $P_1, P_2$, labeled $H_1, H_2$. There is also an outer Hill's region about both $H_1, H_2$, disconnected from them, labeled, $H_O$. $P$ can pass between $H_1$ and $H_2$ when $C_2 < C \lessapprox C_1$.\footnote{$a \lessapprox b$ means $a < b$ and $a$ is slightly less than $b$. Similarly for $a \gtrapprox b$.}  It goes through a small channel region between $H_1, H_2$, labeled, $\mathcal{C}_1$. Similarly, when $C$ reduces further to $C \lessapprox C_2$, $P$ can also pass though a small channel region between $H_2$ and $H_O$, labeled $\mathcal{C}_2$. 

Within the $\mathcal{C}_i$, $i=1,2$,  are unstable retrograde periodic orbits, $\gamma_i$, Lyapunov orbits, whose linearized flow has eigenvalues $\lambda_i = \pm a_i  \pm\sqrt{-1} \ b_i$, where $a_i, b_i, i=1,2$,  are positive real numbers (see \cite{Conley:1968, Conley:1969}).  These periodic orbits have local two-dimensional stable and unstable manifolds, $W^s(\gamma_i), W^u(\gamma_i), i=1,2$, topologically equivalent to cylinders. The trajectories that pass between the Hill's regions move through in interior of the manifolds, as shown in Figure \ref{fig:BigGeom}. The trajectories asymptotically spiral to $\gamma_i$ on $W^s(\gamma_i)$, as $t \rightarrow \infty$, and to $\gamma_i$ on $W^u(\gamma_i)$ as $t \rightarrow -\infty$.
In the case where $C_2 < C <  C_1$, only $\gamma_1, \ W^{s,u}(\gamma_1)$ exist.
 \begin{figure}
\centering
\includegraphics[width=0.65\textwidth, clip, keepaspectratio]{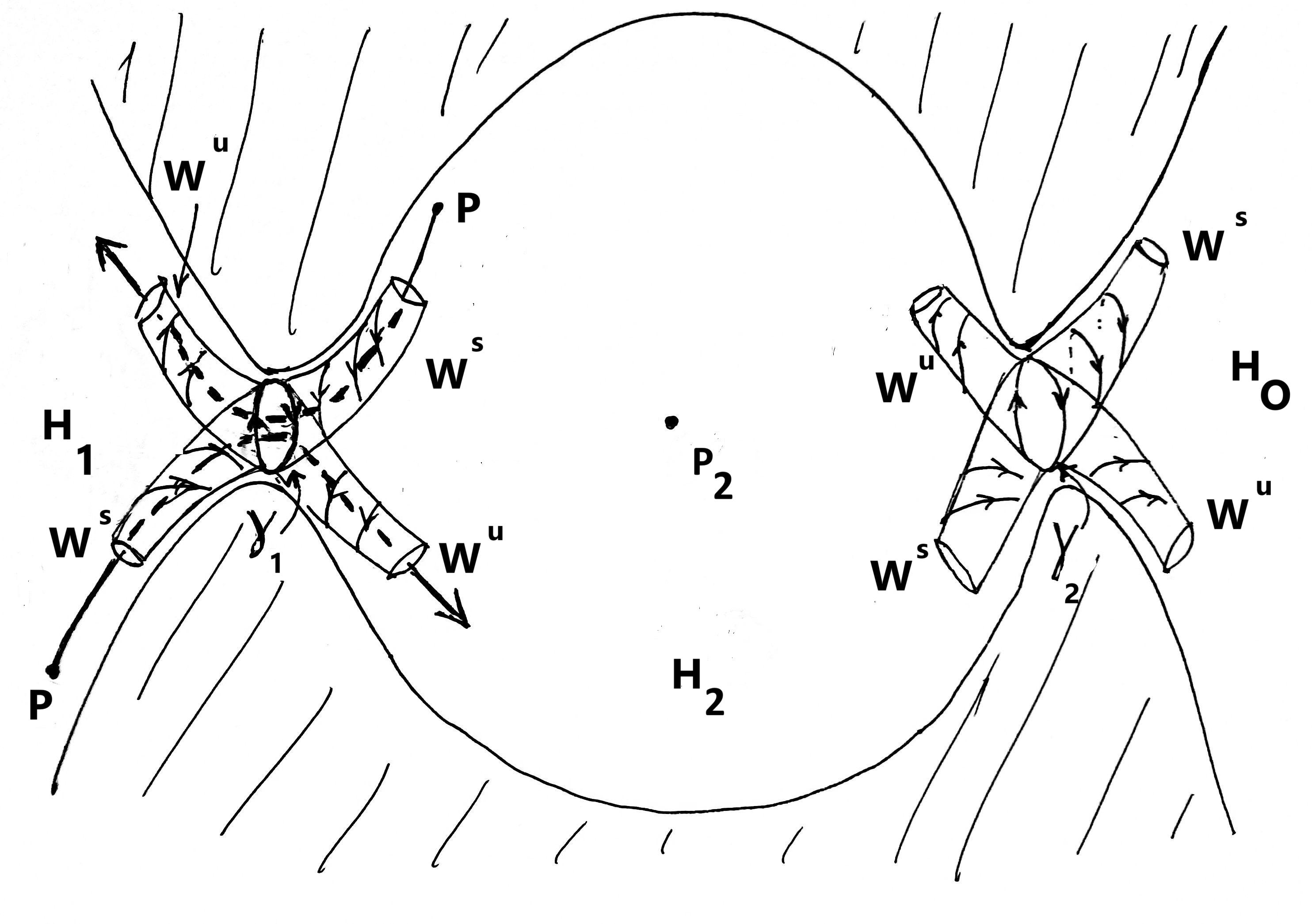}
        		\caption{$W^{s,u}(\gamma_i)$, $i=1, 2$, and the flow of the trajectories towards or away from $\gamma_i$, $C \lessapprox C_2$.  $P$ is shown moving between $H_1, H_2, H_O$ on transit orbits.  This is an illustrative figure.}
	\label{fig:BigGeom}
\end{figure}

\section{Weak Stability Boundary and Related Sets} 
\label{sec:W} 
\medskip\medskip

\subsection{Basic Definitions}
\label{subsec:Definitions}

The definitions of the weak stability boundary are recalled from \cite{BelbrunoGideaTopputo:2010, GarciaGomez:2007, Belbruno:1987}. \\

The center of the $(y_1, y_2)$ coordinate system is translated to $P_2$ by ${\bf y} \rightarrow {\bf Y}: \  Y_1 = y_1 - 1, Y_2 = y_2$, where $Y_1=Y_2=0$.  $P_1$ is located at $Y_1=-1, Y_2 =0$. Let ${\bm{\psi}} (t) = ({\bf{Y}}(t), \dot{\bf{Y}}(t) )$ be a solution of the restricted problem for $P$.   

The Kepler energy between $P$ and $P_2$ is labeled $E_2$. It is a function of  $({\bf{Y}}, \dot{\bf{Y}})$ and given by (\ref{eq:KepEnergyRot}) in Appendix \ref{subsec:KeplerEnergy}.   We define the Kepler energy along the trajectory as $E_2(t) = E_2(\bm{\psi}(t) )$.  

The Jacobi surface, $\sigma$, defined for $P_1$-centered coordinates,  is labeled, $\Sigma = \Sigma(C)= \{({\bf{Y}}, \dot{\bf{Y}}) \in \mathbb{R}^4  | \tilde{J}=C\}$ in $P_2$-centered coordinates, where $\tilde{J}$ is $J$ in $P_2$-centered coordinates given by (\ref{eq:JacobiSCentered}).
\medskip

Consider a line $L(\theta)$ extending from $P_2$ making an angle $\theta$, $0 \leq \theta \leq 2\pi$ with the $Y_1$-axis (see Figure \ref{fig:Cycle}).  Choose initial conditions for ${\bm{\psi}} (t)$ of (\ref{eq:DEs}), in $P_2$-centered coordinates, along each point of $L(\theta)$ at an initial time, $t_0$.  It is assumed that for each initial point on $L(\theta)$, $P$ is at the osculating periapsis of a Keplerian ellipse, where the velocity vector is perpendicular to $L(\theta)$.  

The magnitude $v$ of the velocity ${\dot{\bf{Y}}}$ at $t=t_0$ is determined at each value of $r =  |{\bf{Y}}| > 0$ along $L(\theta)$, so that the initial periapsis value has a given eccentricity $e$, $0 \leq e < 1$. In polar coordinates,
\begin{equation}
v = (\mu(1+e)/r)^{1/2}  - r .
\label{eq:vp}
\end{equation}

Thus, for these initial conditions, $E_2 <0$.  
The corresponding initial osculating periapsis value is $r = a(1-e)$ and $a$ is the initial semimajor axis. The initial value of $E_2 = \mu(e-1)/(2r)$ at $t=t_0$.  Thus, each of these initial conditions along $L(\theta)$ are the periapsis of an osculating ellipse of eccentricity $e$ where $\dot{r} =0$. The initial points lie on the set, 
$$\Lambda = \{ (r, \theta, \dot{r}, \dot{\theta}) \in \mathbb{R}^4| \dot{r} = 0, E_2 < 0 \}.$$  We assume posigrade motion.\footnote{The retrograde case is obtained by the symmetry of solutions.} \\

When $t$ increases from $t_0$, two types of cycling motions are defined: \\

\noindent 
{\it Stable cycling motion} is defined for integer $n$ cycles, $n \geq 1$, when ${\bf{Y}}(t)$ makes $n$ complete cycles about $P_2$, without going around $P_1$, and  $E_2<0$ on $L(\theta)$ on each cycle. It is assumed all intersections of the stable cycling trajectory with $L(\theta)$ are transverse.\footnote{The assumption of transverse intersections for stable cycling motion is generally observed based on numerical evidence (see \cite{BelbrunoGideaTopputo:2010}).} \\


\noindent
{\it Unstable cycling motion} for $n$ cycles occurs when stable cycling for $n$ cycles does not occur.  \\

Unstable cycling for $n$ cycles can occur when (i)  \ ${\bf{Y}}(t)$ makes one full cycle around $P_1$ prior to $n$ full cycles about $P_2$, cycling  around $P_2$ at most $n-1$ times  (see Figure \ref{fig:Cycle}), $n \geq 1$.  Other types of unstable motion about $P_2$ can also occur when: (ii) \ $E_2 \geq 0$ on $L(\theta)$ on any cycle up to the $n$th cycle about $P_2$, (iii) Intersections of the cycling trajectory on $L(\theta)$ are not transverse for all cycles,  (iv) \ $P$ does not return to $L(\theta)$ prior to $n$-cycles and doesn't cycle around $P_1$.  This can occur, for example, for those $n$-unstable trajectories that start on $L(\theta)$ and lie on $W^s(\gamma_i)$; these will neither cycle about $P_2$ or $P_1$. They will asymptotically approach $\gamma_i$. It is numerically observed that (i) is generic. \\

Stable cycling motion for $n$ cycles is denoted by {\it $n$-stable} motion, and unstable cycling motion for $n$ cycles is denoted by {\it $n$-unstable} motion.   1-stable and 1-unstable motion are shown in Figure \ref{fig:Cycle}.  
\medskip

This procedure of determining $n$-stable and $n$-unstable cycling motion defines an numerical algorithm that can be implemented on the computer (see \cite{GarciaGomez:2007, TopputoBelbruno:2009, BelbrunoGideaTopputo:2010}).  This is referred to as the {\it W-algorithm}.\\
\medskip

\begin{figure}
\centering
\includegraphics[width=0.65\textwidth, clip, keepaspectratio]{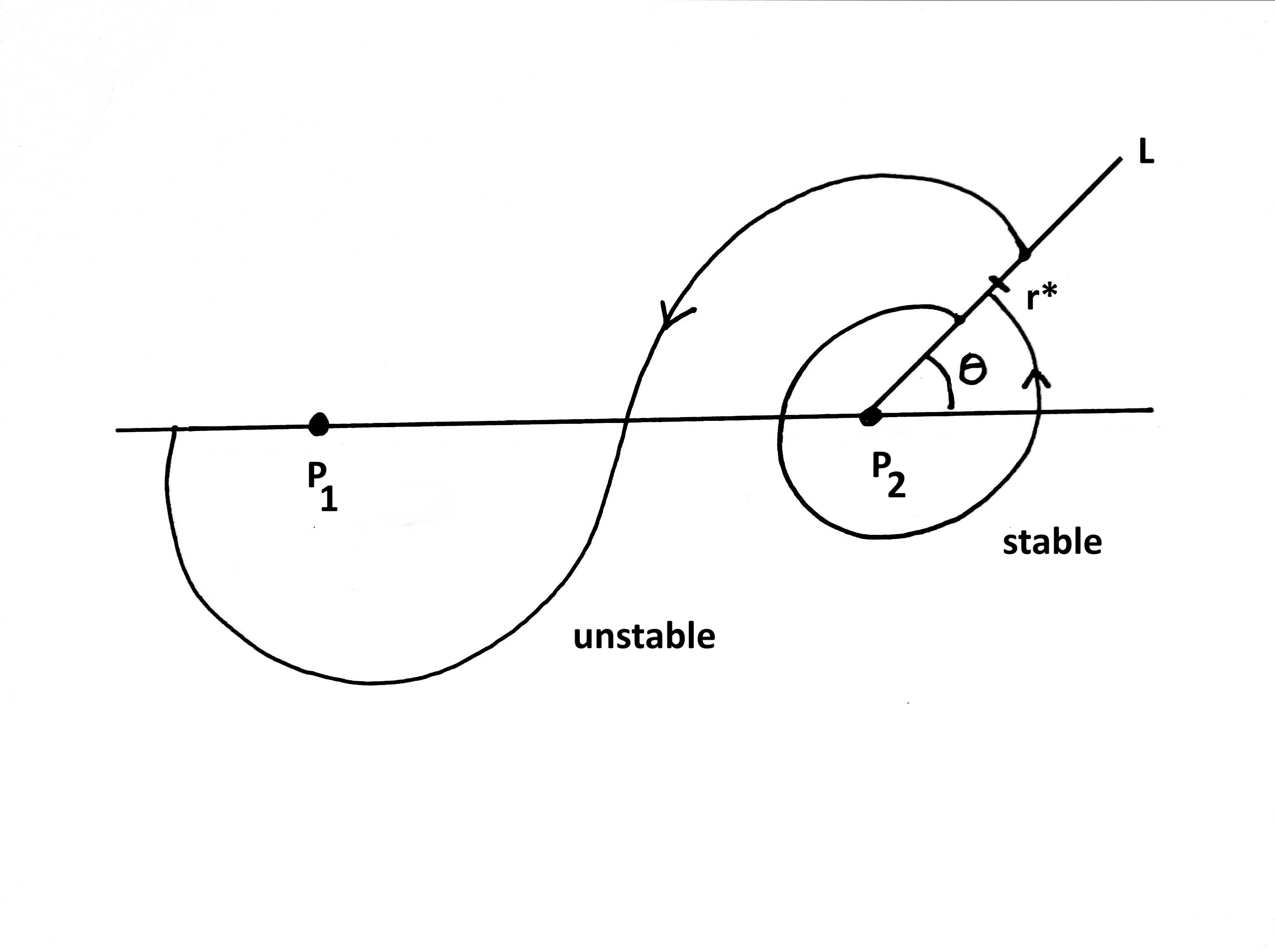}
	\caption{Line $L(\theta)$  emanating from $P_2$ making an angle $\theta$ with the $Y_1$-axis.  1-stable and 1-unstable motion for $n=1$ cycles for $P$ on the trajectory ${\bf{Y}}(t)$.   This is an illustrative figure.}
	\label{fig:Cycle}
\end{figure}

It is noted that in the definition of stable cycling and unstable cycling, 
when initial conditions are chosen for trajectories ${\bm{\psi}} (t)$ along $L(\theta)$ for the W-algorithm, $C$ will vary. For initial conditions sufficiently close to $r=0$, $C$ is arbitrarily large,  and ${\bf{Y}}(t)$ remains in $H_2$ for all $t >t_0$ since $H_1, H_2$ are disjoint.  Thus, in this case, ${\bf{Y}}(t)$ cannot move about $P_1$ and unstable motion by going around $P_1$ is not possible.  \\

It is remarked that the open condition $E_2<0$ is motivated for the sake of applications since such a condition is used for capture of objects (e.g. spacecraft) at $P_2$ \cite{Belbruno:1987}. It serves as a mathematical condition to define stable cycling and as a normalized initial condition for both stable or unstable cycling. The initial condition $\dot{r}=0$ for stable or unstable cycling is used as a consistent normalizing condition, where other values of $\dot{r}$ could have been used.  Using the periapsis of an osculating ellipse is well defined.\\ \\

Sets of all points about $P_2$ leading to $n$-stable and $n$-unstable motion, for a fixed $\theta,e,n$, can be defined. These are called the {\it n-stable set}, {\it n-unstable set}, respectively.   More precisely,
the $n$-stable and $n$-unstable sets on $L(\theta)$ for each $\theta \in [0, 2\pi)$, $e \in [0, 1)$ are denoted by $S_n(\theta, e)$, $U_n(\theta, e)$, respectively.  \\

For a fixed line, $L(\theta)$, and fixed $e, n$, the set $S_n(\theta, e)$ is open, since the $n$-stability condition is open. This is because the condition $E_ 2 < 0$  is an open condition, and  if for some $r>0$ the motion is $n$-stable, then there exists a small $\delta >0$ where it is also $n$-stable in a neighborhood $(r- \delta, r+ \delta)$ along $L(\theta)$ by smooth dependence of solutions of the differential equations on initial conditions.  Since $S_n(\theta, e)$ is an open subset of $L(\theta)$, then topologically this implies that $S_n(\theta, e)$ can be written as a countable union of disjoint open intervals,
\begin{equation}
S_n(\theta,e) =\bigcup_{k \geq 1}  (r^*_{2k-1}, r^*_{2k})  = (r^*_1, r^*_2) \cup  (r^*_3, r^*_4) \cup  \ldots   ,
\label{eq:StableIntervals}
\end{equation}
where $r^*_1 = 0, r > r^*_1$.  The end points of the open intervals $r^*_j$, $j= 1,2, \ldots$, are $n$-unstable. The open intervals in (\ref{eq:StableIntervals}), $(r^*_{2k-1}, r^*_{2k})$,  are referred to as
$I^n_k,  k \geq 1$ for reference.

If a point is $m$-stable then it is $n$-stable for $n < m$, by definition.  Thus, 
\begin{equation}
S_m(\theta, e) \subset S_n(\theta, e), \   n < m.
\label{eq:Sinequal}
\end{equation}
This implies that the more times $P$ cycles about $P_2$ in a stable way, then generally fewer points are available to support that. \\

A change of stability across one of the end points of $(r_{2k-1}^*, r_{2k}^*)$  means that for all $r \in (r_{2k-1}^*, r_{2k}^*)$ the motion is $n$-stable and there exists $r^*, r^{**} \notin  (r_{2k-1}^*, r_{2k}^*)$ arbitrarily close to
to $r_{2k-1}^*$, $r_{2k}^*$, respectively, for which the motion is $n$-unstable.  \\

The $n$-weak stability boundary, $W_n$, is defined to be the union of all $r^*(\theta, e)$ along $L(\theta)$ where there is a change of stability between $n$-stable and $n$-unstable motion,
\begin{equation}
W_n= \{r^*(\theta, e), \ \theta \in [0, 2\pi], e \in [0,1) \} .
\label{eq:Wn}
\end{equation} 
$W_n(\theta,e)$ is a slice of $W_n$. (\ref{eq:Wn}) can equivalently be written as $W_n = \{ W_n(\theta,e), \ \theta \in [0, 2\pi], e \in [0,1) \}$.   Figure \ref{fig:Cycle} shows a point, $r^*$, of $W_1(\theta,e)$.\\

(\ref{eq:Wn}) represents a union of $r^*(\theta, e)$ for the range of $\theta, e$ given. This is equivalent to  
$$W_n = \bigcup_{\theta \in [0, 2\pi], e \in [0, 1)} r^*(\theta, e)$$.

The notation for union in (\ref{eq:Wn}) is used throughout the paper for other sets as well. \\

 $W_n(\theta,e)$ is a discrete set of points for each fixed $\theta, e$ since it's points are the boundary points of disjoint open intervals, $W_n(\theta,e) = \partial S_n(\theta,e)$.  Therefore, taking the union over all $\theta \in [0, 2\pi], e \in [0, 1)$ yields,

\begin{equation}
W_n= \partial S _n. 
\label{eq:WnBdry1}
\end{equation}
$S_n = \{ S_n(\theta,e) ,  \ \theta \in [0, 2\pi], e \in [0,1) \}$ is the general $n$-stable set over all $\theta,e$. \\

It is remarked that (\ref{eq:Sinequal}) does not imply a similar property for $W_n(\theta, e)$.\\


The following sets are defined by taking unions over different parameters, 
 \begin{equation}
S_n(\theta) =  \{ S_n(\theta,e), \   e \in [0, 1) \}, \   S_n(e) =\{ S_n(\theta, e), \ \theta \in [0, 2\pi] \}, 
\label{eq:Sversions}
\end{equation}
\begin{equation}
S(e) = \{ S_n(e),  \  n \geq 1 \}, \ S(\theta) =  \{ S_n(\theta),  \  n \geq 1 \},   \  S = \{S(\theta), \ \theta \in [0, 2\pi]  \}                                
\label{eq:S}
\end{equation}
\medskip

In a similar way, $W_n(\theta), W_n(e), W(e), W(\theta), W$ are defined, as well as for $U_n(\theta, e)$. (It is remarked that also $S = \{S(e), e \in [0,1) \}$.)   \\

$W$ is the {\it general} weak stability boundary and $S$ is the {\it general} stable set. By construction, $C$ varies over $W, S$. Thus these sets are not invariant under the flow of the restricted three-body problem.  $W$ exists over a set of Jacobi energy surfaces, $\Sigma(C)$ in four dimensional phase space.\\

 Note, in $W(\theta)$ the points are obtained as $e$ varies continuously. Therefore, the points of $W(\theta)$ are not discrete, whereas the points of $W(\theta,e)$ are discrete.\\ \\ 
\medskip\medskip

$S_n(\theta,e), W_n(\theta,e)$ have been numerically computed for a fine grid of $\theta \in [0,2\pi]$, for a range of $e$, $n$ and different $\mu$:  \\

\noindent {\it Numerical Results A\\

I. \  It is numerically observed that $S_n(\theta,e), W_n(\theta,e)$ should exist for all finite $n$ and that (\ref{eq:Sinequal}) is also true (see \cite{TopputoBelbruno:2009}, Figure 6, \ $\mu = .00095$; \ $S_n(e)$, $n=1,2, ..., 8; \ e=0$, and Figure 7, \ $W_n(e), n=1; \  e = 0, .2, .3, ..., .95$; \ $n=1,2,..., 8; \ e=0$,  and see \cite{GarciaGomez:2007}, Figures 1,3, $\mu =.01215; \  e=0,.3,.6,.9; \ n = 3,6$).   \\

II. \ It is numerically observed  that $S_n(e)$, is bounded by $W_n(e)$ (see \cite{TopputoBelbruno:2009}, Figure 4, \ $\mu = .00095; \ n=1; \  e=0$; and see  \cite{BelbrunoGideaTopputo:2010}, Figure 4, \ $\mu = .01$; \ $W_1(e), \   e=.02, .6, .95$; \  \cite{JianYiSui:2015}, Figure 4, $\mu = .0000003$; \ $e = .2, .4, .6, .9, .95$).\\

III. \ It is numerically observed that $C$ varies for different points of $W_n, S_n$. For example, the variation of $C$ over $S_1(e)$ for $\mu = .00095$; \ $e=.95$ is shown in (\cite{TopputoBelbruno:2009}, Figure 5). }\\

The observed numerical results in \cite{GarciaGomez:2007, TopputoBelbruno:2009, BelbrunoGideaTopputo:2010, JianYiSui:2015} are for finite values of $n$ not exceeding $n=8$.  Computing these sets is a lengthy numerical process. Visualizing these sets for larger values of $n$ and more $e$ values is not necessary for the purposes of this paper. The general structure and how it changes  as $n$ increases and as $e$ approaches $1$ is generally understood. For example, as $n$ increases, the sets have fewer points. {For the low values of $n$ used in these studies one can get a reasonable idea of the general properties of these sets; however, higher values of $n$ could reveal much more detail and also increasing the refinement of grids used for $e, \theta$ could also be useful, as well as additional values of $\mu$.} \\

{\it  Notation in different papers} - It is important to note that in \cite{TopputoBelbruno:2009} and \cite{BelbrunoGideaTopputo:2010}, the symbol $\mathcal{W}$ is used in place of $S$ and $\mathcal{W}^*$ is used in place of $W$.  In \cite{GarciaGomez:2007}, $\bar{\mathcal{W}}$ is used for $S$ and $\partial{\bar{\mathcal{W}}}$ for $W$.  In \cite{JianYiSui:2015},  $W_n$ is used for $S_n$ and $\partial W_n$ for $W_n$ as in this paper.\\
\medskip\medskip

The sets of most relevance for this paper are $S(\theta), W(\theta)$, for a fixed value of $\theta$, and other variations, such as $S_n(\theta,e), W_n(\theta,e), S_n(e), W_n(e)$ and also $S, W$, which will be clear from the context.\\
\medskip\medskip

\subsection{Sets $W$, $S$, $\hat{S}$, $W'$ and their Properties}
\label{subsec:Properties}
\medskip

Properties of $W, S$ and related sets, $\hat{S}$, $W'$  are described. \\

The following property is true from the basic definition of $W, S$.\\

\begin{prpty}
   $W(\theta)= \partial S(\theta)$.   The trajectories with initial conditions on $L(\theta)$ in $W(\theta)$ are $n$-unstable, $n \geq 1$.
\label{prpty:1}
\end{prpty}

It is numerically observed that  $S_n, W_n$ lie in a region about $P_2$ out roughly to the $L_{1,2}$ distance, then along the approximate distance of $P_2$ to $P_1$, up to $90^o$ ahead and behind $P_2$ (see \cite{TopputoBelbruno:2009, GarciaGomez:2007, BelbrunoGideaTopputo:2010,  JianYiSui:2015}). 

\begin{prpty} 
The trajectories of initial points on $L(\theta)$ in $S(\theta)$ cycle around $P_2$ in an '$n$-stable manner' for $n \geq 1$, i.e. up to the $n$th cycle about $P_2$ the trajectory intersects $L(\theta)$ transversally with $E_2 < 0$.
\label{prpty:2}
\end{prpty}

Let $C_3 \leq C_a < C < C_2$, where $C_a$ is the largest value of $C$ where the Lyapunov orbits, $\gamma_i$, $i=1,2$, don't collide with $P_2$ and the channels 
$\mathcal{C}_i$ exist to allow $P$ to escape into $H_1$ or $H_O$.   For example, in \cite{BelbrunoGideaTopputo:2010}, where $\mu = .01215$ the value of $C_a = 3.15$ is assumed, where, in that case, $C_1 \approx 3.20034, C_2 \approx 3.18416, C_3 \approx 3.02415$.

\begin{prpty}
$S(\theta)$ is open and  bounded on $L(\theta)$ for $C_a < C < C_1$.  
\label{prpty:3}
\end{prpty}

\noindent
 Proof  -
\medskip

\noindent
Case 1: \   $C_2 \leq C < C_1$ .  $H_2$ has one channel  $\mathcal{C}_1$ about $L_1$ connecting $H_2$ with $H_1$ of width $\mathcal{O}(\mu^{1/3})$ for $\mu$ sufficiently small \cite{Conley:1968, LlibreMartinezSimo:1985}. This enables unstable cycling about $P_2$ where $P$ can escape $H_2$ and eventually cycle about $P_1$. 
$S_n(\theta)$ is bounded. This follows by first noting that $L(\theta) \cap H_2$ is bounded for all $\theta \in [0, 2\pi]$. Thus for any trajectory starting on $L(\theta)$ on $\Lambda$, it will always be bounded, which includes $n$-stable cycles.  For $L(\theta) \cap H_O$, all trajectories with initial conditions on $\Lambda$ that cycle about $P_2$ are $n$-unstable. This follows since for any cycle about $P_2$, $P$ automatically makes a full cycle about $P_1$ when it returns to $L(\theta)$.\footnote{If $C_a$ were allowed to go below $C_3$, $P$ need not make a full cycle about $P_1$.} This implies $P$ cannot make $n$-stable cycles about $P_2$ in this case. Thus, $S_n(\theta)$ is bounded for all $n$, and so is $S(\theta)$. This set is also open since it is the infinite union over $n \geq 1$, $e \in [0, 1)$ of open sets, $S_n(\theta,e)$, and therefore open. \\

\noindent
Case 2: \  $C_a \leq C < C_2$ .    The only difference in this case is that another channel opens up, $\mathcal{C}_2$, connecting $H_2$ to $H_O$. This proof that $S(\theta)$ is bounded and open follows in a similar way to Case 1. \\

  This proves Property \ref{prpty:3}.  {This proof is rigorous, not requiring special assumptions or numerical results.}
\medskip\medskip

 It is noted that numerical results {verify} $S$ is bounded from figures for the stable sets for different $n, \  e$ (see, Numerical Results A, and,  for example, \cite{TopputoBelbruno:2009, GarciaGomez:2007}).
\medskip\medskip

   Also, numerically, it would be difficult to discern $S_n(\theta,e)$ on $L(\theta)$ for large $n$. This follows since $S_n(\theta,e)$ can be written as a countable union of open intervals and because of the property (\ref{eq:Sinequal}) of $S_n(\theta, e)$, the $r^*$ values would be difficult to numerically discern. 
This situation becomes more complicated for different values of $e \in [0,1)$. Since the solutions of the differential equations vary continuously as a function of the initial conditions on $L(\theta)$, then near the intervals for a fixed $e = e_0$, as $e$ varies  $(e_0 -\delta, e_0 + \delta)$, $\delta$ small,  the open intervals defining  $S(\theta,e)$ vary by a small amount. Thus, the $r^*(e)$ vary by a small amount. This yields overlapping open intervals and boundary points for $e$ near $e_0$, which would be difficult to discern. \\
\medskip \medskip\medskip\medskip

\noindent {\it Infinite Cycling Stable Motion, $\hat{S}(\theta)$}
\medskip\medskip

It is of interest in this study to find a subset of points of $S_n(\theta,e)$ on $L(\theta)$ where the trajectories cycle about $P_2$ in a stable manner for all time, i.e. for all infinitely many cycles. This is defined by, 
\begin{equation}
\hat{S}(\theta, e) =  \lim_{n \rightarrow \infty} S_n(\theta,e) \ = \bigcap_{n \geq 1}  S_n(\theta, e)  .
\label{eq:SLimit}
\end{equation}
{\it This is a set of all  points starting on $L(\theta) \cap \Lambda$ where trajectories through these points cycle infinitely often about $P_2$ in a stable manner for all time. $C$ varies over this set. $\hat{S}(\theta) = \{ \hat{S}(\theta, e), \  e\in [0, 1) \}$. }  

\begin{prpty}
$\hat{S}(\theta, e)$  need not be open. Assuming it is non-empty, then it is bounded for $C_a < C < C_1$.  
\label{prpty:4}
\end{prpty}

\noindent
Proof - 
\medskip

 $S_n(\theta, e) $ form a  sequence of open sets satisfying  (\ref{eq:Sinequal}) that are descreasing in size, $n =1,2, \ldots$.  As $n \rightarrow \infty$, the limiting
intersection need not be open since the intersection of infinitiely many open sets need not be open.  

 The boundedness follows from Property \ref{prpty:3}.  This proves Property \ref{prpty:4}. {This proof is rigorous, not requiring numerical results}.
\medskip\medskip

If $C$ is sufficiently large, then it can be analytically proven $\hat{S}(\theta)$ is non-empty.

\begin{thm}
If $C$ is sufficiently large ($ C \in C_L $) then $\hat{S}(\theta)$ is non-empty. 
\label{thm:KAM}
\end{thm}

\noindent 
Proof - 
\medskip

The existence of two-dimensional KAM tori about $P_2$ on each three-dimensional energy surface $\Sigma(C)$ is proven in \cite{Kummer:1979} for $C$ sufficiently large, in a set $C_L = \{C \ | \ C > M \}$ for $M$ sufficiently large.  For $C \in C_L$, $H_2$ shrinks down to $r \gtrapprox 0$. Say,  $0 < r < \delta(M)$, $\delta$ is small.  

The tori separate the phase space on $\Sigma(C)$. Trajectories on the tori, $\bm{\psi}(t) = ( {\bf Y}(t), \dot{ \bf{Y}}(t) )$, with initial conditions at say $t=t_0$, have two frequencies of motion.  One is rotational frequency $\omega = 1$, and the other is the osculating frequency $\omega^*(t, \mu)$  along the trajectory. To first order, $\omega^*(t,0)= \tilde{\omega}$ is a constant, representing the elliptic Kepler frequency of motion for an ellipse  of constant semimajor axis $a \gtrapprox 0$ and constant eccentricity $e \in [0,1)$.  $\omega / \tilde{\omega}$ are sufficiently irrational satisfying diophantine conditions. The motion between the tori is not sufficiently irrational and chaotic. 

Since the tori separate phase space, the trajectory $\bm{\psi}(t)$ for $P$ is trapped between the tori for all time, where $a,e$ vary by $\mathcal{O}(\mu)$. $a,e$ also vary by $\mathcal{O}(\mu)$ when $\bm{\psi}(t)$ lies on the invariant tori. Thus, in the full three-dimensional region on $\Sigma(C)$ where the tori exist near $P_2$, the periapsis and apoapsis distances of $P$ from $P_2$ vary as $a(1-e) + \mathcal{O}(\mu), a(1+e) + \mathcal{O}(\mu)$, respectively, for all time. In particular, $E_2(\bm{\psi}(t)) < 0$.   The family of tori transversally cut a surface of section $S_{\theta}$ (see  (\ref{eq:SectionDef})), in a family of invariant curves.  

The line $L(\theta)$ extending from $P_2$  is given by the $r$-axis of the section. The trajectories transversally cut $L(\theta)$ as they cycle about $P_2$ on or between the KAM tori.

This structure implies that trajectories moving about $P_2$ for all time for $C \in C_L$ starting on $\Lambda$ will lie on or between the KAM tori. In particular, trajectories with initial conditions on $\hat{S}(\theta)$ will move on or betwen KAM tori. Thus, $\hat{S}(\theta)$ is non-empty.\\

 This proves Theorem \ref{thm:KAM}. {As can be seen, this proof is carried out rigorously, not based on numerical work or any assumptions .}
\medskip\medskip\medskip

It is noted that numerical results show that $\{ S_n(\theta, e) \ | \ n \geq 1 \}$  is non-empty for the ranges of $C$ given in Property \ref{prpty:4}.  (See Numerical Results A).  Although the numerical investigations only go as far as $n =8$,  it is seen that the figures of $S_n(\theta, e)$ appear consistent as $n$ varies, and points of $S_n(\theta, e)$ are close to $P_2$. These suggest that $S_n(\theta, e)$ is nonempty as $ n \rightarrow \infty$, although this is not proven.  It is assumed that $\hat{S}(\theta,e)$ is nonempty for this paper. It is out of the scope of this paper to numerically investigate the structure of $\hat{S}(\theta,e)$ for this range of $C$, and is an interesting topic to study further.\\ 

{The assumption that $\hat{S}(\theta,e)$ is nonempty is validated by the results in \cite{KoonLoMarsdenRoss:2000}. They show the existence of resonant periodic orbits about $P_2$ for $\mu = .00095$ for the range of $C$ of relevance in this paper.  This orbit therefore cycles infinitely often about $P_2$ in a stable manner. This orbit for $P$ is for the comet Oterma modeled in the planar circular restricted three-body problem used in this paper. The proof in \cite{KoonLoMarsdenRoss:2000} is semi-analytic; that is, based on numerical results, together with theory.}    \\

When $C \notin C_L$, the results from \cite{Kummer:1979} on the use of the KAM theorem cannot be used to study bounded motion for all time from initial values on $\hat{S}$.
However, interesting properties of  $\hat{S}$ can be obtained for $C_a < C < C_1$ without the use of KAM tori.  This described in Section \ref{subsec:GlobalS}. 
\medskip\medskip\medskip\medskip

\noindent {\it Notation} \ \ Since $C$ varies over  $\hat{S}(\theta)$, then points of $\hat{S}(\theta)$ on $\tilde{J}=C$ can be denoted as, $\hat{S}_C(\theta) = \hat{S}(\theta) \cap \Sigma(C)$.  Similarly, $W_C(\theta)  = W(\theta) \cap \Sigma(C), \ S_C(\theta) = S(\theta) \cap \Sigma(C)$, etc.  The subscript $C$ can be omitted if the context is clear.\\ \\

\noindent {\it Properties of $W(\theta)$ and Infinite Cycling Unstable Motion, $W'(\theta)$}
\medskip\medskip\medskip

The set  $W_n(\theta,e)$ is considered.\\  

\begin{prpty}
$W_n(\theta, e)$ has a countable set of points.  $W(\theta, e)$ contains a countably infinite number of points,  and  $W(\theta)$ has an uncountably infinite number of points.  $\dim W_n(\theta, e) = \dim W(\theta,e) = 0$
\label{prpty:5}
\end{prpty}

\noindent
Proof -  
\medskip

$W_n(\theta, e)$ contains a countable set of points, as follows from (\ref{eq:StableIntervals}).  Thus, $W(\theta,e)$ contains a countably infinite number of points since it is the union over $n \geq 1$.  

Taking the union over all $e$, $W(\theta) = \{ W(\theta, e), \   e \in [0, 1) \}$, implies that $W(\theta)$ contains an uncountably infinite number of points.

The dimension of $W_n(\theta,e)$ is zero since it contains a countable set of points. Taking the union over all $n = 1,2, \ldots$, $W(\theta, e)$ also has zero dimension. This proves Property \ref{prpty:5}. {This proof is analytic, not requiring numerical results or assumptions.}\\  \\

It is remarked that  $\dim W(\theta)$ is not clear at this stage.  Although $\dim W_n(\theta,e) = 0$, when the union is taken over all $e \in [0,1)$ and $C$ is restricted to a given range, further analyisis is needed.
\medskip

It is noted that $W_n(\theta,e)$ may contain a countably infinite number of points. This can happen when the lengths of the open sets, $I^n_j, j=1,2, \ldots$ comprising $S_n(\theta,e)$ converge to $0$. \\

\begin{lem}
 $W(\theta)$  is bounded for $C_a < C < C_1$. $W(\theta)$ is closed assuming uniform convergence of all sequences $r^*_k(\theta,e), \ e \in [0,1), k = 1,2, \ldots.$
\label{lem:ClosedBdd}
\end{lem}

 \noindent
Proof -   
\medskip

$W(\theta)$ is bounded since $S(\theta)$ is bounded by Property \ref{prpty:3} and $W(\theta) = \partial S(\theta)$ by Property \ref{prpty:1}.  

Since $W_n(\theta,e)$ is a countable discrete set, it is closed since its boundary itself.  This is also true for $W(\theta,e)$.  It is not automatically the case that
$W(\theta)$ is closed by taking the infinite union over $e \in [0,1)$ of closed sets,  since that need not be closed. 
 
$W(\theta)$  is closed if it contains all of its limit points.   This is shown to be the case by considering the closed set $W(\theta,e)$. Since this set is closed, any convergent sequence of points $\{ r^*_k(\theta,e) \}$ in $W(\theta,e)$ converges to a limit point $\alpha(\theta,e) \in W(\theta, e)$, as $k \rightarrow \infty$. 

 More precisely, given an $\epsilon > 0$ there exists a $K(\theta, e, \epsilon) > 0$, bounded for each $e \in [0,1)$,  such that for $k > K$,  $|r^*_k(\theta,e) - \alpha(\theta,e)| < \epsilon$, where $\alpha(\theta,e) \in W(\theta,e)$.
 This is a pointwise convergence  for each $e \in [0,1)$. For $W(\theta)$ to be closed, this convergence should be uniform in $e$. That is,   $K$ is independent of $e$.  This is achieved by setting $\tilde{K}(\theta, \epsilon) = \sup_{e \in [0,1)}K(\theta, e, \epsilon)$ and assuming $K(\theta, e, \epsilon)$ is uniformly bounded for all $e \in [0,1)$.
 
 Thus, any convergent sequence of points $\{ r^*_k(\theta) \}$ in $W(\theta)$ converges to a limit point $\alpha(\theta) \in W(\theta)$ as $k \rightarrow \infty$.  This implies $W(\theta)$  is closed.
 \\

 This proves Lemma \ref{lem:ClosedBdd}. \\

{The assumption required in the proof of this lemma is, 
\begin{assum}
The sequences $ r^*_k(\theta), \ k= 1,2, \ldots$ converge uniformly for $e \in [0,1)$. 
\label{assum:1}
\end{assum}
}

This assumption is suggested from numerical studies {for low values of $n$ up to $8$, where computing $W_n(\theta)$ is straight forward for $e = .95$ \cite{TopputoBelbruno:2009}} but proving this numerically requires computing $W_n(\theta,e)$ for high values of $n$ for $e  \uparrow  1$, which is out of the scope of this study. {It is possible that an analytic proof could be done, but none currently exists.} 

{ If this assumption were not true, then $W(\theta)$ would not be a closed set. The validity of Assumption \ref{assum:1} and hence the validity of Lemma \ref{lem:ClosedBdd} enable the proof of Property \ref{prpty:7}, in the following:}
 \\

In an analogous manner to the definition of $\hat{S}(\theta,e)$, points common with $W_n(\theta,  e)$ for all $n \geq 1$ can be defined as, 
\begin{equation}
W'(\theta, e) = \lim_{n \rightarrow \infty} W_n(\theta,e)  =   \bigcap_{n \geq 1} W_n(\theta,  e) 
\label{eq:WPrime}
\end{equation}

It is remarked that points in $W'(\theta, e)$ give rise to trajectories, ${\bm \psi}(t)$, that cycle about $P_2$ infinitely often as $ t \rightarrow \infty$. Since these initial points yield unstable cycling motion, 
this means that, by definition,  in the process of cylng about $P_2$, ${\bm \psi}(t)$ can also cycle about $P_1$, before returning to $L(\theta)$, or $E_2 \geq 0$ on at least one intersection with $L(\theta)$. The possibility of $P$ not returning to $L(\theta)$ cannot occur in this case, whether cycling around $P_1$ or not. \\

\begin{prpty}
$W'(\theta)  = \partial \hat{S}(\theta)$.
$W'(\theta, e)$ is countable. $W'(\theta)$ is bounded.  $W'(\theta)$ is closed, assuming uniform convergence of all sequences $r^*_k(\theta,e), \ e \in [0,1), k = 1,2, \ldots$.  $C_a < C < C_1$. 
\label{prpty:7}
\end{prpty}

\noindent
Proof -  
\medskip

$W'(\theta) = \partial \hat{S}(\theta)$: \ \  Consider the line $L(\theta)$ for some $\theta \in [0, 2\pi]$ and fixed $e\in [0,1)$.   
For a fixed $n \geq 1$, $W_n(\theta,e)= \partial S_n(\theta,e)$ (see (\ref{eq:WnBdry1})).  Taking the limit as $n \rightarrow \infty$,  yields
\begin{equation}
W'(\theta,e) = \partial \hat{S}(\theta,e)
\label{eq:ResultingLimit}
\end{equation}
as follows from  (\ref{eq:SLimit}), (\ref{eq:WPrime}).   (The boundary points of $\hat{S}(\theta,e)$, that comprise $W'(\theta,e)$,  represent boundary points of the remaining open intervals, $I^n_j, j \geq 1$, as $n \rightarrow \infty$ of $S_n(\theta,e)$. These points are the boundary points of the initial points for trajectories that cycle about $P_2$ for all time in a stable manner.)

(\ref{eq:ResultingLimit}) is true for each $e \in [0,1)$. Thus, $W'(\theta) = \partial \hat{S}(\theta)$.   

$W'(\theta, e)$  is a countable set by Property \ref{prpty:5} since $W_n(\theta,  e)$ is countable and the countable intersection of countable sets is countable.

$W_n(\theta,e)$ is closed and bounded by Lemma \ref{lem:ClosedBdd}. Thus, $W'(\theta, e)$ is closed  since (\ref{eq:WPrime}) is a countable intersection of closed sets which is therefore closed.  Property \ref{prpty:4} implies it is bounded.  This is true for each $e \in [0,1)$.  Thus, $W'(\theta,e)$ is closed and bounded for each $\theta \in [0, 2\pi], e \in [0,1)$.  $W'(\theta)$ is closed following a similar argument as in the proof of Lemma \ref{lem:ClosedBdd} for $W(\theta)$, under the same assumption, and the boundedness of $W(\theta)$ implies the same for $W'(\theta)$.\\ 

This proves Property \ref{prpty:7}.
\medskip

{As is seen, Assumption \ref{assum:1} enables the proof of Property \ref{prpty:7}. This plays a key role in this paper. It is used to show that $W'(\theta)$ is a Cantor set for each value the Jacobi constant, $C$, in Lemma \ref{lem:CantorSet}, which is a key result of this paper.  If this assumption were not true then the Cantor structure of $W'$ would not be necessarily true, but in that case some restricted variant of it would likely be true. For example, the range of $e$ could be restricted so that it is bounded away from $e=1$ and $W'$ would have a partial Cantor structure.}
\medskip\medskip\medskip

It is noted that $\dim W'(\theta)$ is not clear.  This is the case, since from Property \ref{prpty:5}, $\dim W_n( \theta,e) = 0$, and (\ref{eq:WPrime}) implies $\dim  W'(\theta,e) = 0$.   However, when the union is taken over $e \in [0, 1]$ and for $C_a < C < C_1$,  the dimension of $W'(\theta)$ requires further analysis.  The dimension of $W'(\theta)$ is determined in Section \ref{subsec:GlobalW}.\\

Similar to $\hat{S}(\theta)$, $W'(\theta)$ is defined over a range of $C$.  It is noted that  the value of $C$ of a point $w \in W'(\theta)$ is in general different than the values of $C$ for points in $\hat{S}(\theta)$ adjacent to it that it bounds.      \\

\noindent Definition \ \  The set of points of $W'(\theta)$ that lie on the energy surface $\Sigma(C)$ are given by $W'_C(\theta) = W'(\theta) \cap \Sigma(C)$.  
\\

It is noted that $\hat{S} \subset \Lambda, \ W'(\theta) \subset \Lambda$.\\

\noindent
{\it Dynamical Summary A \  Trajectories, $\bm{\psi}(t)$,  with initial conditions at $t=0$  on  $\hat{S}(\theta), W'(\theta)$ satisfy,\\

a.)  $\hat{S}(\theta)$ defines initial points at $t=0$ on $L(\theta) \cap \Lambda$ for trajectories cycling $P_2$ infinitely often, for all $t>0$, in a 'stable manner' ($E_2 <0$ on all intersections with $L(\theta)$ and are transversal where $\dot{r}$ need not be zero). \\

b.)  $W'(\theta)$ defines points on $L(\theta) \cap \Lambda$ that bound $\hat{S}(\theta)$, and initial points for trajectories cycling $P_2$ for all $t>0$ in an 'unstable manner' (For example, $E_2 \geq 0$ on an intersection with $L(\theta)$), or intersections with $L(\theta)$ can be nontransversal, or the trajectory cycles about $P_1$ before returning to $L(\theta)$.) } 
\medskip\medskip

It is noted that the open sets $I^n_k$ comprising $S_n(\theta,e)$ change for different $n$ and $e$. It was briefly noted in \cite{GarciaGomez:2007} that $S_n(\theta,e)$ suggestive of a Cantor set.
It turns out that the points of $S_n(\theta,e)$ are not in a Cantor set, but rather the boundary points $W'(\theta)$ of $\hat{S}(\theta)$ . Topologically,  these are the boundary points of $\lim_{n \rightarrow \infty} S_n(\theta)$.  This is proven in Section \ref{subsec:GlobalW} using a main result in \cite{BelbrunoGideaTopputo:2010} on relating $W_n(\theta)$ to the invariant manifolds $W^{s}(\gamma_i)\ i=1,2$. \\ \\

\subsection{$\bm{W_n(\theta)}$, $\bm{W'(\theta)}$, Invariant Manifolds, Cantor Sets}
\label{subsec:GlobalW}

The relationship of $W_n(\theta)$ to invariant manifolds is described.
\medskip\medskip

Three cases for the range of the Jacobi integral are considered. \\ 

Case 1:  $C_1 \leq C$, \ \  Case 2:  $C_2 \leq C < C_1$, \ \  Case 3:  $C_a < C < C_2$.  It is assumed that $\mu > 0$ is small, $\mu \gtrapprox 0$.\\

A Poincar\'e surface of section $S_{\theta}$ of two dimensions can be defined on $\Sigma(C) \cap H_2$ for the flow of the differential equations in Cases 1,2,3.
$S_{\theta}$  is defined on the surface $\Sigma(C)$ for each fixed $\theta \in [0, 2\pi]$ and for each fixed $C$.  The section has coordinates $(r, \dot{r})$, where $\dot{\theta} > 0$. It is given by, 
\begin{equation}
S_{\theta} = \{(r, \dot{r})| \  \theta = \theta_0, \ \dot{\theta} = \dot{\theta}(r, \dot{r}, \theta)> 0, r> 0\} 
\label{eq:SectionDef}
\end{equation} 
(see Appendix \ref{subsec:SolvingforThetaDot}).\footnote{$S_{\theta}$ not to be confused with $S(\theta)$. }  \\

\noindent
{\it Notation} \   $S_{\theta}$ is used for a range of $C$ values or a single $C$ value which will be clear from the context. Sometimes the notation, $S_{\theta}(C)$, is used for a particular $C$ value.  \\

The flow of the differential equations defines a two-dimensional Poincar\'e map, $\Phi$, on $S_{\theta} \cap \Sigma(C)$ for each $C$ in the given range, and each $\mu \gtrapprox 0$,
$$ \Phi :  S_{\theta}\rightarrow S_{\theta} .$$

Case 1 - $H_2$ lies within the distance to $L_1, L_2$.  $\gamma_i, \  i=1,2$ do not exist, and $\mathcal{C}_i$ do not exist for each fixed energy surface $\Sigma(C)$ in the $C$ range. 
The range of $C \in C_L$ in Theorem \ref{thm:KAM} is a subset of this case.\\

Case 2 - $P$ can escape $H_2$ and move between $H_1, H_2$ through $\mathcal{C}_1$. For each fixed $C$ in the range,  the invariant manifolds $W^{s,u}(\gamma_1)$ lie on $\Sigma(C)$ and they extend into $H_2$.

 The first intersection (cut) of $W^s(\gamma_1)$ and $W^u(\gamma_1)$ on $S_{\theta}(C)$ are topological circles.  If the manifolds have a tranvserse intersection, then the topological circle breaks up on the section (see \cite{GideaMasdemont:2007}) under the flow of the differential equations.  It is assumed the trajectories satisfy topological properties (see Section \ref{subsec:TrajProp}) and for details, see \cite{BelbrunoGideaTopputo:2010}. 


After infinitiely many intersections of $W^s(\gamma_1)$ and $W^u(\gamma_1)$ on $S_{\theta}(C)$,  a hyperbolic invariant set is obtained resulting from the intersections. The existence of the hyperbolic invariant set on $S_{\theta}(C)$ follows by the Smale-Birkhoff theorem (see \cite{Moser:1973, GuckenheimerHolmes:1983}). This hyperbolic set is a Cantor set of hyperbolic points.
This set is labeled, $\mathcal{C}_C(\mu)$, or just $\mathcal{C}_C$ for brevity.  If the union is taken over a range of $C$ in $[C_2, C_1)$ and $\mu \gtrapprox 0$, in the $(C, \mu)$-plane where the manifolds have transverse intersection, then the symbol $\mathcal{C}$ is used. This represents a union of Cantor sets lying on a set of energy levels $\Sigma(C)$, for  $(C,  \mu)$. \\ 

Let $p \in \mathcal{C}_C$.  $p$ has one-dimensional stable and unstable manifolds on $S_{\theta}(C)$.
These manifolds intersect at all the hyperbolic points.  This creates a complicated network of hyperbolic points and manifolds (tangle).   $\mathcal{C}_C$ is closed, self-similar, totally disconnected, and  has zero topological dimension on $\Sigma(C)$, for each $(C, \mu)$ where transversality occurs.  \\

  $\Phi: \mathcal{C}_C \rightarrow \mathcal{C}_C$  is chaotic on $\mathcal{C}_C$ for each $C$ with a dense set of periodic points ($q$  a periodic point implies
 $\Phi^m(q) = q$, for some $m>0$).  There also exists a dense orbit, $\{\Phi^k(\tilde{p}), \ k \geq 0\}$ for some $\tilde{p}$ \ ($\Phi^0(p) = p$, \ $\Phi^2(p) = \Phi \circ \Phi (p)$).   
\\

What can be said about the occurance of transversality of $W^u(\gamma_1), W^s(\gamma_1)$ in $H_2$ in this case?

An analytic proof (see \cite{LlibreMartinezSimo:1985})  answers this in the $H_1$ region, where it is proven that for a discrete set of $C \lessapprox C_1, \ \mu \gtrapprox 0$,  $\{ (C_k, \ \mu_k), \ k=1,2,  \ldots \}$,  $W^u(\gamma_1), W^s(\gamma_1)$ indeed have transverse homoclinic intersection.   $(C_k, \mu_k)$ are numerically determined for a set of $k$ in \cite{LlibreMartinezSimo:1985}.\\

It is reasonable to assume that these manifolds have transversal intersection in $H_2$.  In fact, as is described in Case 3, transversal heteroclinic intersections  of invariant manifolds are proven to exist when both channels exist about $L_1, L_2$, which is a more general situation than considered in Case 2.\\

Case 3 - When $C$ is less than $C_2$, $P$ can move through $\mathcal{C}_1, \mathcal{C}_2$. Also, $\gamma_1, \gamma_2$ grow in size.

Instead of considering a transverse homoclinic intersection of $W^{s}(\gamma_1), W^{u}(\gamma_1)$ in $H_2$, it can be proven that there is a transverse heteroclinic intersection between
$W^s(\gamma_1)$ and $W^u(\gamma_2)$ in $H_2$, which by symmetry of solutions under the map, $Y_2 \rightarrow -Y_2, \ \dot{Y}_1 \rightarrow -\dot{Y}_1$, also gives transverse intersection between $W^u(\gamma_1)$ and $W^s(\gamma_2)$.  

This heteroclinic intersection is proven to exist, in a semi-analytic proof (numerically assisted) in the interesting paper by J. Marsden, et. al. \cite{KoonLoMarsdenRoss:2000}.  It is done in the case of $P_1=$Sun, $P_2=$Jupiter, where $\mu = .00095$, and for specific $C \lessapprox C_2$ ($C=3.037$, \  $C_2 =3.03836$). It is shown that $W^s(\gamma_1)$ and $W^u(\gamma_2)$ transversally intersect on a surface of section.  This choice of $C, \mu$ is motivated by the trajectory for the resonance transitioning comet Oterma studied in \cite{KoonLoMarsdenRoss:2000}.   As noted in \cite{KoonLoMarsdenRoss:2000}, their analysis is valid for other resonant transition comets about the Sun, for the same $\mu$, but different $C \lessapprox C_2$.\footnote{A list of many other similar resonance comets about the Sun in resonance with Jupiter is given in \cite{BelbrunoMarsden:1997, Belbruno:2007}.} This yields other values of $C$, for the same $\mu$, where transversality of  $W^s(\gamma_1), W^u(\gamma_1)$ would similarly exist.  More generally, it is stated in \cite{KoonLoMarsdenRoss:2000} that it was numerically observed that transverse intersection of the manifolds occurs for many values $C_a < C < C_2, \mu \gtrapprox 0$. \\

Thus, it can be assumed, based on semi-analytic analysis, that heteroclinic intersections exist for set of $C_a < C < C_2 , \mu \gtrapprox 0$ in the $(C,\mu)$-plane.  

     As in Case 2,  there exists an invariant Cantor set for $\Phi$, also labeled $\mathcal{C}_C$, with the properties given in Case 2.  $\mathcal{C}_C$ exists on each three-dimensional energy surface, $\Sigma(C)$ for set of $C_a<C < C_2, \mu \gtrapprox 0$.  \\ 

{It is noted that the following assumption is used in Cases 2, 3, }
{
\begin{assum}
The manifolds intersect transversally.
\label{assum:2}
\end{assum}
}

{The transversality of the manifolds $W^s(\gamma_1), W^u(\gamma_2)$  in Case 3 and $W^s(\gamma_1), W^u(\gamma_1)$ in Case 2 is necessary for the use of the Smale-Birkhoff theorem to deduce the existence of Cantor sets of hyperbolic points. Without this transversality, the Cantor structure of $W'$ cannot be deduced. } 
\medskip

{ It is seen that the proof of the transversality of the intersection of the manifolds of $W^s(\gamma_1)$ and $W^u(\gamma_2)$ for a set of $C$ within the interval $(C_a, C_2)$ {for Case 3} follows from a semi-analytic proof in \cite{KoonLoMarsdenRoss:2000} for $\mu =.00095$.  This semi-analytic proof can be extended for many more values of $\mu$ and $C$ as stated in \cite{KoonLoMarsdenRoss:2000}.  This proof is not purely analytic, but does guarantee transversality for some restricted cases of $C$, $\mu$. In Case 2, where $C_2 < C < C_1$, a proof would be needed in that range, but from \cite{LlibreMartinezSimo:1985} transversality of $W^s(\gamma_1)$ and $W^u(\gamma_1)$ should be true.  Even though Case 3 has a semi-analytic proof for many cases of $C$, $\mu$ it does not explore all cases of small $\mu$ nor all $C$ in the given range. An analytic proof for transversality of the manifolds is not currently available that covers all the relevant $\mu$, $C$ ranges. To make the proof of transversaility of the manifolds in both cases rigorous, an analytic proof would be required, which is an open problem. It is reasonable to conjecture that an analytic proof is possible.  } 
\medskip\medskip

Definition \ \  Let $I_C \subset I = \{ C \ | \ C_a < C < C_1 \}$ and $I_{\mu} \subset  I_{\delta} = \{\mu \ | \ 0 < \mu < \delta \}$, where $\delta \ll 1$ is sufficiently small,  be sets in the $(C, \mu)$-plane, where, within $C_2 \leq C  < C_1$, 
$W^s(\gamma_1) , W^u(\gamma_1)$ in $H_2$ have transverse homoclinic intersection; and for those values within  $C_a < C < C_2$, $W^s(\gamma_1), \ W^u(\gamma_2)$ in $H_2$ have transverse heteroclinic intersection. This set of $(C,\mu)$ values is labeled $V_{C,\mu} =  \{ C \in I_C, \mu \in I_{\mu} \}$.
\\ \\

\noindent {\it Relationship of $W_n$ to $W^s(\gamma_i), i=1,2$ }  \\

  The main result of \cite{BelbrunoGideaTopputo:2010} is that $W_n$ lies on $W^s(\gamma_i), i=1,2$, under some conditions.  Assumptions are necessary on the trajectories space, described in Appendix \ref{subsec:TrajProp}, referred to as Hypothesis A.  The main technical result is recalled, \\
\medskip\medskip

\noindent
{\it Main result in \cite{BelbrunoGideaTopputo:2010} }
\medskip\medskip

\noindent 
{\it Let $W_n^A$ be the union of the $n-1$ stable manifold intersections with $S_{\theta}$ on $\Lambda$, for all  $\theta \in [0, 2\pi]$ for each $C \in I, \mu \in I_{\delta}$, then taking the union over all $C \in I$ for each fixed $\mu \in I_{\delta}$, where it assumed the trajectories satisfy Hypothesis A.  Then, $W_n^A = W_n$. }\\
\medskip

\

{
It is noted that the following assumption is required in the previous result stated more precisely in (\ref{eq:WnSet}),  
\begin{assum}
The trajectories satisfy Hypothesis A in Appendix \ref{subsec:TrajProp}.  
\label{assum:3}
\end{assum}
}

{This assumption is made for the following reason.  When $W^s(\gamma_i), i=1,2$ extend from the Lyapunov orbits, $\gamma_i$, respectively, in $H_2$, they are observed to wind around $P_2$ for the assumed values of $C$.  When the $W$-algorithm is applied along a line $L(\theta)$, it intersects the manifolds, and the $n$-unstable points turn out to correspond to the interior points of the stable manifolds on $L(\theta)$ while the $n$-stable points correspond to points in $L(\theta)$ exterior to the manifolds.The boundary between these sets of points on $L(\theta)$  correspond to the manifold points. By restricting to trajectories satisfying Assumption \ref{assum:3}, the $W$-algorithm is well defined along $L(\theta)$. This assumption also rules out points where the $W$-algorithm is not well defined. This leaves gaps in $W_n(\theta)$.  This assumption works fine for the lower values of $n$ used in several papers  that goes up to $n=8$ (see Numerical Results A, Section 3.1). Based on this it is logical to assume it should be valid for any $n$ which is the assumption. If this were not true, then (\ref{eq:WnSet}) may not be valid for all $n$.}

{The previous result is} stated more precisely {as}, \\

\begin{equation}
W_n =  W_n^A  = \{ (r, \dot{r}, \theta, \dot{\theta}) \in \{ W^s_{\theta, n-1}(\gamma_1) \cup W^s_{\theta, n-1}(\gamma_2) \}  \cap \Lambda, \ C \in I,   \  \theta \in [0, 2\pi]\},
\label{eq:WnSet}
\end{equation}
for each fixed $\mu \in I_{\delta}$, where $W^s_{\theta, n-1}(\gamma_i)$ is the $(n-1)$st intersection of $W^s(\gamma_i)$ with $S_{\theta}$. \   Fixing $\theta$ defines $W_n(\theta)$. Assumption \ref{assum:3} is satisfied. \\

\medskip

Fixing $\theta$,  (\ref{eq:WnSet}) implies,

\begin{lem}
\label{lem:KeyRelation}
\begin{equation} 
W_n(\theta) =  W_n^A(\theta). 
\label{eq:KeyEqu}
\end{equation}
\end{lem}
\noindent for $C \in I$ for each fixed $\mu \in I_{\delta}$.
\medskip

This equates the $W_n(\theta)$ points from the $W$-algorithm with the points $W_n^A$ of the stable manifolds. 

Since $W(\theta) = \{ W_n(\theta) , \ n \geq 1\}$, and $W^A(\theta) = \{  W_n^A(\theta) , n \geq 1\}$ then (\ref{eq:KeyEqu}) implies \ \  $W(\theta)  = W^A(\theta)$.
\medskip

\begin{cor} 
\label{cor:1}  Assume $C \in I_C$, then for each $\mu \in I_{\mu}$, 
\begin{equation}
 W'(\theta) =  \lim_{n \rightarrow \infty} W^A_n(\theta) = \mathcal{C} \cap \Lambda
\label{eq:KeyEqu2}
\end{equation}
on $S_{\theta}$.  $\mathcal{C}$ is a union for Cantor sets for
each $\mu$.
\end{cor}

Proof -\ \ (\ref{eq:KeyEqu}) relates $n$ cycles of the manifolds $W^s(\gamma_i), \ i=1,2$, about $P_2$, intersecting $S_{\theta}$, on the right hand side of (\ref{eq:KeyEqu}),  to  $W_n(\theta)$ for $n$ cycles of trajectories about $P_2$ starting on $\Lambda$ on the left hand side of (\ref{eq:KeyEqu}).  This is true for each $\mu \in I_{\delta}$ and for all $C \in I$.
\medskip

\noindent   
Letting $n \rightarrow \infty$ for infinite cycling, the left hand side of  (\ref{eq:KeyEqu}) yields $W'(\theta)$  and the right hand side of (\ref{eq:KeyEqu}) yields the limit set $\mathcal{C}$ on $\Lambda \cap S_{\theta}$ assuming $C \in I_C \subset I, \ \mu \in I_{\mu} \subset I_{\delta}$.  That is, for $(C, \mu) \in V_{C,\mu}$.  

This proves Corollary \ref{cor:1}.\\   

{This corollary is proven rigorously, but it assumes previous results that used assumptions already discussed. Without those previous results this would not be true.}
\medskip

 It is remarked that by Corollary \ref{cor:1}, for a fixed $C$ and $\mu$, $W_C'(\theta)= \mathcal{C}_C \cap \Lambda \cap S_{\theta}(C) $. \\

Since $\mathcal{C}_C$ is a Cantor set, and $\mathcal{C}_C \cap \Lambda \subset \mathcal{C}_C$, then $\dim W'_C(\theta) =0$. This answers the question on the dimension of $W'(\theta)$ after the proof of Property \ref{prpty:7}.\\

Numerical results demonstrate that $W_n^A(\theta)$ is non-empty \cite{BelbrunoGideaTopputo:2010}.  Thus, by Corollary \ref{cor:1}, $W'(\theta)$ is assumed to be non-empty.   It would be interesting to numerically visualize the  Cantor set structure of $W'(\theta)$. This is beyond the scope of this paper. \\

{It is noted that the demonstration $W_n^A(\theta)$ is non-empty is numerically based and done for a limited number of values of $n$, $C$, $\mu$.  In this sense, it isn't completely rigorous.} \\

\noindent
Definition  \ \  $W'(\theta)$ is the weak stability boundary for infinitely many cycles of $P$ about $P_2$.}  \\

\noindent {\it Topological picture}\\

The Cantor set $\mathcal{C}_C$ exists on $\Sigma(C)$ for each $C$ in $I_C$ and for each $\mu \in I_{\mu}$. It exists on the two-dimensional section $S_{\theta}$, where both $r, \dot{r}$ vary. The section lies in the three-dimensional space $\Sigma(C)$. $\mathcal{C}_C$ is zero-dimensional.   Taking the union for $C \in I_C$, yields a set of Jacobi energy surfaces and a set $\mathcal{C}$ of Cantor sets, belonging to $S_{\theta} \cap \Sigma$, $\Sigma = \{\Sigma(C), C \in I_C \}$.

$W_C'(\theta) \subset \mathcal{C}_C$, \ $W'(\theta)\subset \mathcal{C}$.    The points of $W_C'(\theta)$ and $W'(\theta)$ lie on the $r$-axis of $S_{\theta}$, with $E_2 < 0$. $W_C'(\theta)$ is zero-dimensional.  \\

The next result shows that $W_C'(\theta)$ is a Cantor set, 
   \medskip

\begin{lem}
$W_C'(\theta)$ is a Cantor set for each $\mu \in I_{\mu}$.  $W'(\theta)$ is a union of Cantor sets over $C \in I_C$.  $W' = \{ W'(\theta), \theta \in [0, 2\pi] \}$ is an infinite union of Cantor sets. \\
\label{lem:CantorSet}
\end{lem}

\noindent 
Proof- 
\medskip

  From the remark following Corollary \ref{cor:1}, $W'_C(\theta)$ is a subset of a Cantor set. $W_C'(\theta)$ is closed by Property \ref{prpty:7}. This implies $W'_C(\theta)$ is a Cantor set.  Thus $W'(\theta)$ is a union of Cantor sets. $W'$ is an infinite union of Cantor sets taken over $\theta \in [0, 2\pi]$. 
\medskip 

{As noted in the Introduction, the proof of this lemma uses Assumptions \ref{assum:1}, \ref{assum:2}, \ref{assum:3} and therefore isn't completely rigorous.}
\medskip

\begin{lem}
$W'_C(\theta)$ is not invariant under $\Phi$ .
\label{lem:1}
\end{lem}

This is proven by first noting that since
$W'_C(\theta) \subset \mathcal{C}_C$, it consists of hyperbolic points on $S_{\theta}(C)$ for $\Phi$. These points lie on the $r$-axis, where $E_2 < 0$. 
Although $\mathcal{C}_C$ is an invariant set for $\Phi$,  $W_C'(\theta)$ need not be an invariant set for $\Phi$. 
This is the case since if $p \in W_C'(\theta)$, then the orbit of iterates, $\{ \Phi^k(p), \   k \geq 0 \}$ may remain in $W'_C(\theta)$ for a finite number of iterations, where $\dot{r}=0, E_2 <0$ and then leave $W'_C(\theta)$ with $\dot{r} \neq 0$ or $E_2 \geq 0$ and lie in $\mathcal{C} - W'_C(\theta)$. 
\medskip\medskip\medskip

Assumption \ \ It is assumed $\mu \in I_{\mu}$ is fixed for the remainder of this paper, unless otherwise indicated.\\

\subsection{An Invariant Set for $\bm{\Phi}$}
\label{subsec:GlobalS}
In this section an invariant set for $\Phi$ is obtained by considering $\hat{S}(\theta)$.\\

$W'(\theta)$ and $\hat{S}(\theta)$  exist on the $r$-axis of $S_{\theta} \cap \Lambda$, for $C \in I_C$ and $C \in \hat{I}_C \subset I$, respectively.   $\hat{S}(\theta)$ is  bounded by Cantor points of $W'(\theta)$, where $C$ varies throughout the points of $\hat{S}(\theta)$ and $W'(\theta)$.  Thus, the points of these sets lie on different energy surfaces $\Sigma(C)$.  

Since the points of $W'(\theta)$ are the boundaries of points in $\hat{S}(\theta)$, then there exists values of $C$ within $\hat{I}_C$ that will lie arbitrarily close to values of $C$ in $I_C$. It would be interesting to numerically explore the distribution of points of $\hat{S}(\theta)$ and $W'(\theta)$ on the $r$-axis of the section. This is beyond the scope of this paper.

It is recalled that if a trajectory ${\bm \psi}(t)$ has an initial condition on $\hat{S}(\theta)$  (where $\dot{r}=0, E_2 <0$)  then it cycles around $P_2$ for all $t > 0$. On each cycle, it transversally intersects $S_{\theta}$ with $E_2 < 0$, but not necessarily $\dot{r}=0$.  This implies,

\begin{lem}
$\hat{S}(\theta)$ is not invariant under $\Phi$.
\label{lem:2}
\end{lem}
\medskip\medskip

It is noted that ${\bm \psi}(t)$ has two frequencies of the motion when cycling about $P_2$.  One is $\omega=1$ for the rotating frame. The other is $\tilde{\omega}(t)$ along the trajectory.  
$\tilde{\omega} (t)$ can have significant variation as a function of $t$. \\
\medskip

Definition \ \ The set of trajectories for initial conditions in $\hat{S}_C(\theta)$ move in a bounded region $T_C$ on $\Sigma(C)$ about $P_2$, for the given values of $C$.  
$T = \{T_C, \ C \in \hat{I}_C \}$ is the union of the $T_C$.
\medskip\medskip

Although $\hat{S}(\theta)$ isn't invariant under $\Phi$, an invariant map can be constructed for a set related to  $\hat{S}(\theta)$.  
\medskip\medskip\medskip

\noindent
{\it An Invariant Set $S^*(\theta)$  for $\Phi$}
\medskip\medskip

The set is defined, 
\begin{equation}
S_C^*(\theta) = \left\{ \{\Phi^k(p), k \geq 0\}, \ p \in  \hat{S}_C(\theta) \right\},
\label{eq:SStar}
\end{equation}
$\hat{S}_C(\theta) \subset S_C^*(\theta)$.   
Each point $p \in S_C^*(\theta)$ lies on $\Sigma(C) \cap S_{\theta}$. $S^*(\theta) = \{ S^*_C(\theta), \ C \in \hat{I}_C \}$.
\medskip\medskip

 For a point  $p \in \hat{S}(\theta)$, $\dot{r}=0, E_2 < 0$  and for all subsequent iterates of $p$, $E_2 < 0$ since the trajectory has stable motion relative to $L(\theta)$. \\

Thus, 
\medskip\medskip

\begin{lem}
$ S^*(\theta)$ is invariant under $\Phi$,
\begin{equation}
\Phi : S^*(\theta)  \rightarrow S^*(\theta),
\label{eq:Sstar}
\end{equation}
defined on $S_{\theta}$ for each $C \in \hat{I}_C$, on $\Sigma(C)$. $ S^*(\theta)$ consists of points $p = (r,\dot{r})$ with $E_2 < 0$ on $S_{\theta}$.  $\hat{S}(\theta)$ consists of those points of $S^*(\theta)$ with $\dot{r} =0$. $T$ is  formed from all the trajectories starting in $S^*(\theta)$.
\label{lem:3}
\end{lem}

It is noted that $W'(\theta)$ bounds $\hat{S}(\theta)$ but  $\mathcal{C}(\theta)$ does not necessarily bound  $S^*(\theta)$.    The points of  $W'(\theta)$  bound the points of $\hat{S}(\theta)$ since the points of $W_n(\theta)$ bound the points of $S_n(\theta)$. The points of both of these sets satisfy $\dot{r}=0, E_2 <0$ and lie on the $r$-axis of $S_{\theta}$.  This is not satisfied in general for $S^*(\theta)$. The points of $S^*(\theta), \ \mathcal{C}$  are generally not on the $r$-axis of $S_{\theta}$. Thus, the bounding of  $S^*(\theta)$ by $\mathcal{C}$ doesn't necessarily follow. The relationship between $S^*(\theta)$ and $\mathcal{C}(\theta)$ is not clear and whose study is beyond the scope of this paper.\\

\section{Mappings and Properties}
\label{sec:Geom}

The following assumptions are made,  \\

\noindent
Assumptions A \ \  The energy surface $\Sigma(C)$ is considered and a section $S_{\theta}$,  $S_{\theta} \cap \Sigma(C)$, $C \in I$  is fixed, $\theta \in [0, 2\pi]$ is fixed, $\mu \in I_{\mu}$ is fixed.

\begin{lem}
$\Phi$ is real analytic on $S_{\theta}$. In particular it is real analytic on  $\mathcal{C}, S^*(\theta)$, $ W'(\theta), \hat{S}(\theta)$, 
\label{lem:4}
\end{lem}

\noindent
Proof -
\medskip

 Under Assumptions A, $\Phi$ is a real analytic map at each point of $S_{\theta}$ and in particular at each point of $\mathcal{C} \subset S_{\theta}$, $C \in I_C$, and therefore on $W'(\theta) \subset  \mathcal{C}$.  This follows from the real analyticity of the 
solutions $\bm{\psi}(t)$ for the differential equations as a function of initial conditions on $S_{\theta}$ and $t$.   Similarly, $\Phi$ is real analytic at each point of $S^*(\theta)$, $C \in \hat{I}_C$, and therefore at each point of $\hat{S}(\theta) \subset S^*(\theta)$.\\

\noindent
{\it Mapping of  $\mathcal{C}(\theta), S^*(\theta)$ }
\medskip\medskip

\begin{figure}
\centering
\includegraphics[width=0.85\textwidth, clip, keepaspectratio]{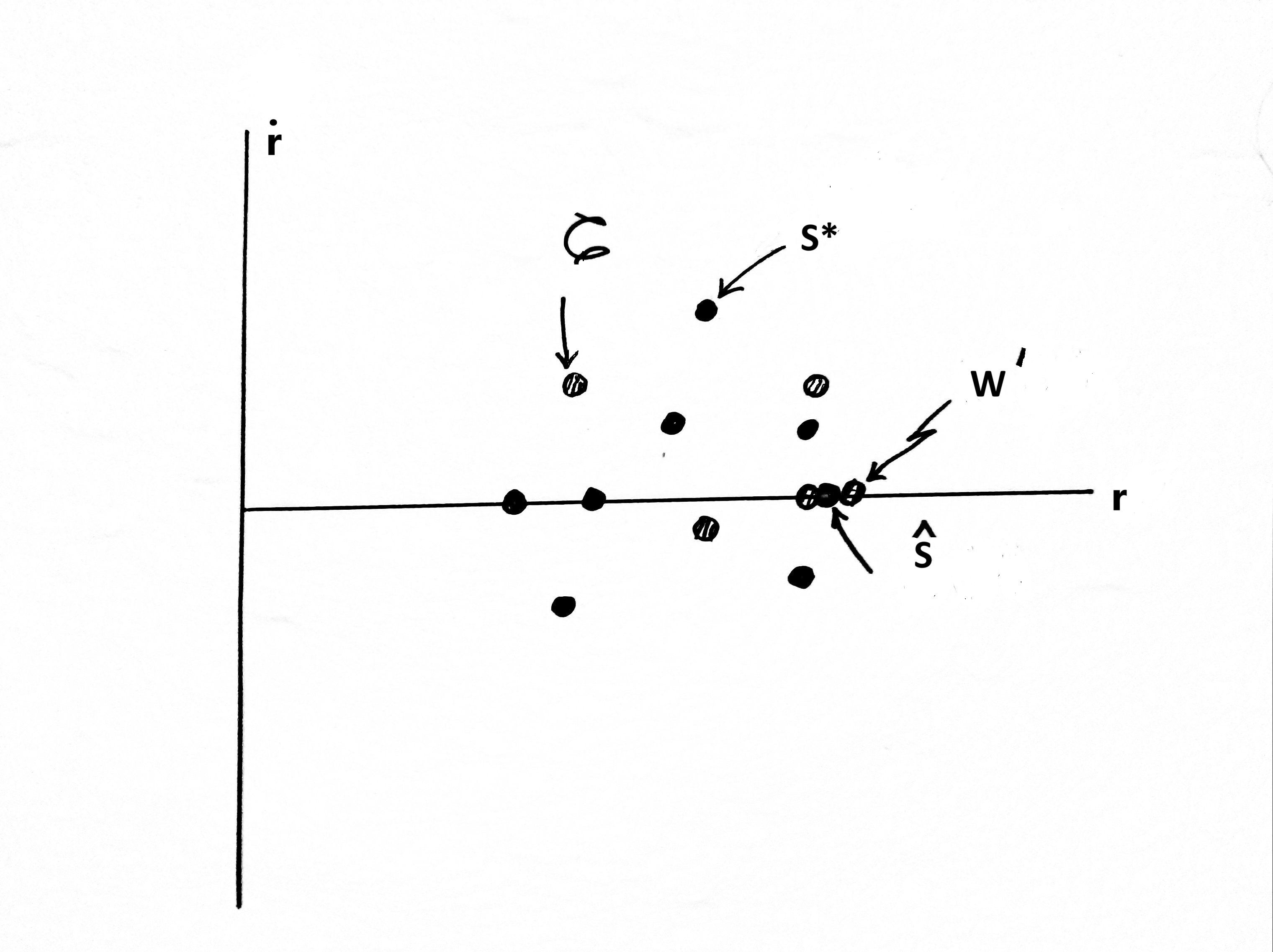}
	\caption{ $\hat{S}(\theta), \  S^*(\theta), \   W'(\theta), \   \mathcal{C}.$  {This is a rough sketch.}} 
	\label{fig:MapGeom}
\end{figure}

$S^*(\theta)$ is invariant under $\Phi$ on $S_{\theta} \cap \Sigma(C)$,  $C \in \hat{I}_C$. $\mathcal{C}$, is invariant under $\Phi$, on $S_{\theta} \cap \Sigma(C)$, $C \in I_C$ (see Figure \ref{fig:MapGeom}).   {This figure is a rough sketch and not generated numerically. There are no numerical simulations of these points, and this is for future work. }  Trajectories with initial values in $S^*(\theta)$ cycle about $P_2$ for all time in a stable manner. By Lemma \ref{lem:3}, this set of trajectories is given by $T$, $S^*(\theta)  =  T \cap S_{\theta}$.    The trajectories of $T$ lie in an three-dimensional annular region about $P_2$ on $\Sigma(C)$.  Trajectories that have initial points  in $\mathcal{C}$ are dynamically unstable since they have two-dimensional stable and unstable manifolds.\\
\medskip  \medskip

\noindent
{\it Mapping of  $W'(\theta), \hat{S}(\theta)$ }
\medskip\medskip

The points of $W'(\theta), \hat{S}(\theta)$ lie on $S_{\theta} \cap \Lambda$, or equivalently on $L(\theta)$ for all $C \in I_C, C \in \hat{I}_C$, respectivly. 

For each fixed $C$,  the iterates of $\Phi$ on $W'(\theta), \hat{S}(\theta)$,  are mapped to points of $\mathcal{C}(\theta), S^*(\theta)$, respectively,  as seen in Figure \ref{fig:MapGeom}.
\medskip\medskip

\medskip

Thus, in summary,\\
\medskip

\noindent $\mathcal{C}, S^*(\theta)$ are invariant under $\Phi$ and $W'(\theta) \subset \mathcal{C}, \  \hat{S}(\theta) \subset  S^*(\theta)$ are not invariant under $\Phi$. $C \in I_C$ for points in $\mathcal{C}$, and $C \in \hat{I}_C$ for points in $S^*$.  $\mu \in I_{\mu}$ is fixed.\\

\noindent
\section {Bounded Regions Defined by $\bm{W'}, \bm{\hat{S}}$ } 
\label{sec:Mandelbrot}
\medskip

$W'(\theta)$, $\hat{S}(\theta)$ are previously described on $S_{\theta} \cap \Lambda$ using both manifolds and the W-algorithm. These sets equivalently exist on the line $L(\theta)$ satisfying $\dot{r}=0, E_2 < 0$. That is, they can be described on the $Y_1,Y_2$-plane for the given value of $\theta$ along $L(\theta)$.  As previously described one gets the limit of the sets $W_n(\theta) = \{ W_n(\theta,e), \ e \in [0, 1) \}$ and  $S_n(\theta) = \{S_n(\theta,e), \ e \in [0, 1) \}$  for $n \rightarrow \infty$ as $W'(\theta) , \ \hat{S}(\theta)$, respectively.  Taking the union over $\theta \in [0, 2\pi]$, yields $W', \ \hat{S}$, respectively.  

\begin{figure}
\centering
\includegraphics[width=0.85\textwidth, clip, keepaspectratio]{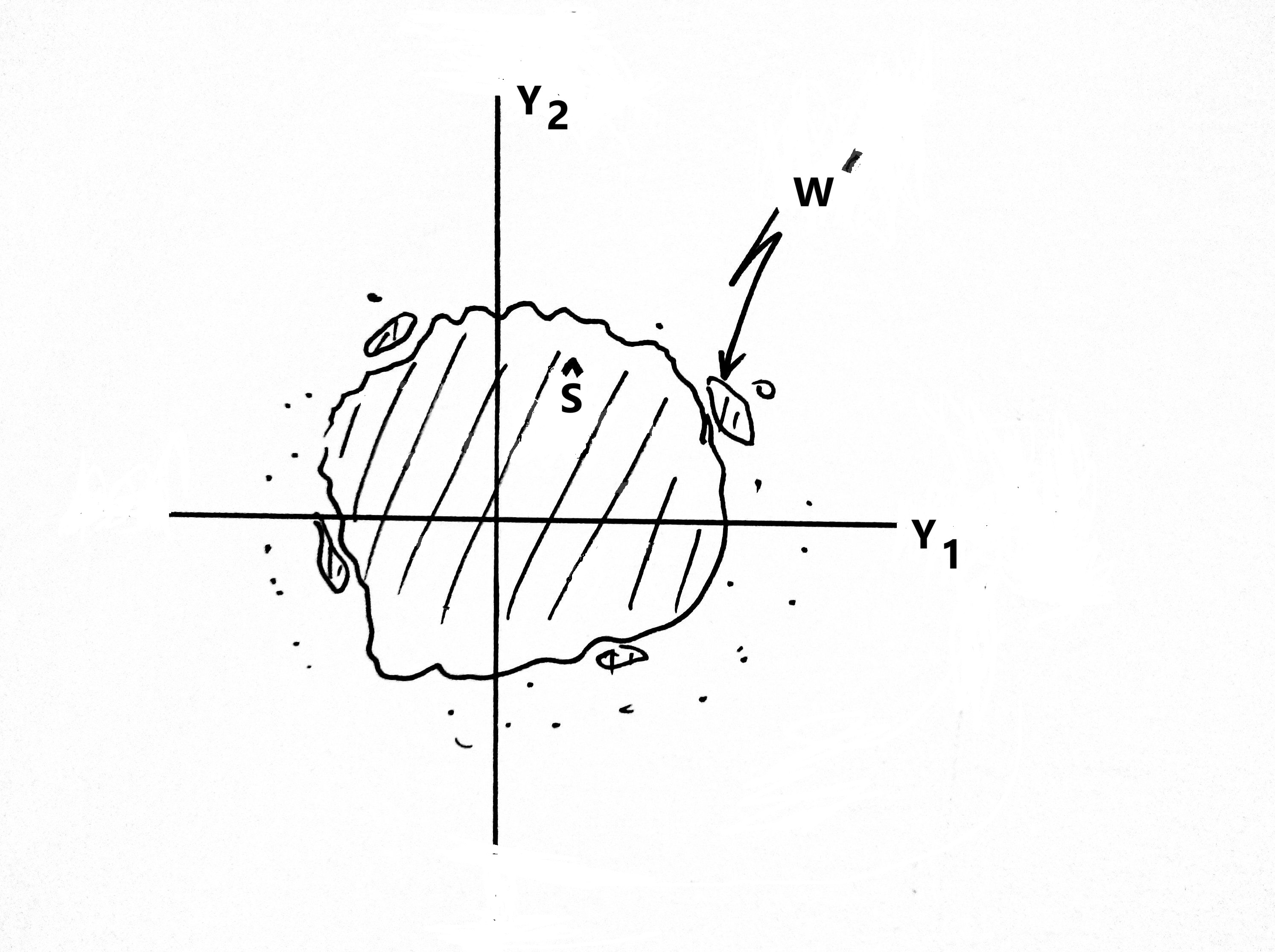}
	\caption{$W', \ \hat{S}$ on $P_{Y, \Lambda}$. 
                {This is a rough sketch.}    }
	\label{fig:WPrimeSHat}
\end{figure}

$W_n(e), S_n(e)$ are described in Numerical Results A for various $n$, $e$.  They show a complex pattern.  One sees this clearly in \cite{TopputoBelbruno:2009} (Figure 4 for $n=1$, $e=0$). It can be seen that $S_1(0)$ fills in a large region about $P_2$, and many smaller regions, whose boundary is $W_1(0)$. Many of the smaller regions are isolated islands. As indicated by other figures in \cite{TopputoBelbruno:2009}, described in Numerical Results A, as $n$ increases the sets decrease on size, seen clearly in \cite{TopputoBelbruno:2009} where $1 \leq n \leq 8$. Different values of $e$ show a similar structure. This is also observed in \cite{GarciaGomez:2007, JianYiSui:2015}.  These figures suggest that as $n \rightarrow \infty$, $W', \hat{S}$ decrease in size significantly for $C \in I_C, \ C \in \hat{I}_C$, respectively.  For a fixed $\theta$, $W'(\theta)$  is a union of Cantor sets, on each respective $\Sigma(C)$, that would be hard to discern numerically, as would the set $\hat{S}(\theta)$, especially if the different $C$ were close in value.
See Figure \ref{fig:WPrimeSHat} which is just a rough sketch (not a numerical simulation) for all $\theta \in [0, 2\pi]$. 

  A numerical investigation of $W', \hat{S}$ for large values of $n$ is beyond the scope of this paper. \\    
\medskip\medskip

\noindent{\it  A Parametric Plane and an Iteration Plane}\\

The points of $W', \ \hat{S}$ can be viewed in the  ${\bf Y} = (Y_1,Y_2)$-plane. Each of these sets are defined by first varying a position in the  $(Y_1,Y_2)$-plane along $L(\theta)$, for a value of $\theta$, then using the W-algorithm to see which points exists along $L(\theta)$ after infinitely many cycles from $W_n(\theta)$ and $S_n(\theta)$.  Since $\dot{r}=0$, then the other velocity $\dot{\theta}$ component (in polar coordinates) needs to be suitably adjusted satisfying $E_2({\bf Y}, \dot{\bf{Y}}) <0$ that yields a desired value of $e \in [0,1)$ (see (\ref{eq:KepEnergyRot})).

Each point of $W', \ \hat{S}$ in the $(Y_1,Y_2)$-plane, $P_{Y_1,Y_2}$,  has associated to it a value of $E_2 <0$ and satisfies, $\dot{r}=0$.  Thus, the points belong to the projection of $\Lambda$ onto $P_{Y_1,Y_2}$. $\Lambda = \{(Y_1,Y_2, \dot{Y}_1, \dot{Y}_2) \in \mathbb{R}^4 \ | \dot{r}({\bf{Y}}, \dot{\bf{Y}})=0, E_2({\bf{Y}}, \dot{\bf{Y}} ) < 0 \}$.\\

\noindent Definition \ \ $P_{Y, \Lambda}$ denotes the points of the $(Y_1,Y_2)$-plane as a projection from $\Lambda$. \\

Thus,  $W'$, $\hat{S}$  exist on $P_{Y, \Lambda}$. $C \in I$ varies at different points of  $P_{Y, \Lambda}$.\\

It is noted that the points of $P_{Y, \Lambda}$ can be viewed as the  periapsis points of Keplerian elliptical orbits.

A key observation is that the values of $W'$, $\hat{S}$ along $L(\theta)$ in $P_{Y, \Lambda}$ are used as points to iterate on $S_{\theta}$, starting on the $r$-axis. 

If one chooses $p \in W'(\theta)$, then $\Phi^n(p) \in \mathcal{C}$ for all $n \geq 1$, where $\dot{r}$ need not be $0$, which lie on $\Sigma(C)$ for the corresponding value of $C \in I_C$. Similarly, if $p \in \hat{S}(\theta)$ then $\Phi^n(p) \in S^*(\theta)$ for all $n \geq 1$, where $\dot{r}$ need not be $0$, which lie on $\Sigma(C)$ for the corresponding value of $C \in \hat{I}_C$.  \\

{\it Thus, 
\medskip

\noindent
(a.) \ The points of $P_{Y, \Lambda}$ provide the initial points for $\Phi$ to iterate on $S_{\theta}$ for a given $\theta$.  $W', \ \hat{S}$ exist on $P_{Y, \Lambda}$.\\

\noindent (b.) \ The $(r, \dot{r})$-plane defining  $S_{\theta}$ is the plane of iterates for $\Phi$ that belong to $\mathcal{C}, \  S^*(\theta)$, respectively (see Figures \ref{fig:MapGeom},  \ref{fig:WPrimeSHat}).} \\ \\ 

\noindent Definition \ \  $P_{Y, \Lambda}$ is referred to as the {\it initial parameter plane} for $\Phi$ and $S_{\theta}$ is referred to as the {\it iteration plane} for $\Phi$.\\

There are points in $P_{Y, \Lambda}$  that don't belong to $W', \ \hat{S}$, and are complementary to these sets on $P_{Y, \Lambda}$.

\begin{lem}
The complementary points to those of $W'(\theta) \cup \hat{S}(\theta)$ relative to $L(\theta)$, on $P_{Y, \Lambda}$, are the unstable points, $U_n(\theta)$, for any finite $n$,  that don't belong to $W'(\theta)$, labeled, $\bar{U}_n(\theta)  = \{p \in (U_n(\theta) - W'(\theta)) \cap (L(\theta) \cap \Lambda), \ n \geq 1 \}$, $C \in I$.
\label{lem:6}
\end{lem}

The proof follows by the definition of $W'(\theta), \hat{S}(\theta)$.
\medskip\medskip

Geometrically the points of $\bar{U}_n  = \{   \bar{U}_n(\theta)  , \  \theta \in [0, 2\pi] \}$ represent all the points of $P_{Y, \Lambda}$ in Figure \ref{fig:WPrimeSHat} outside the hatched regions defining points of $\hat{S}$.  {This figure is a rough sketch and not generated numerically. There are no numerical simulations of these points, and this is for future work. } \\

Trajectories having initial values in $\bar{U}_n$ can eventually move out of $H_2$ into $H_1$ or $H_O$, prior to $n$ cycles of a given trajectory $\bm{\psi}$ about $P_2$ and cycle about $P_1$ and not returning to $H_2$,  or asymptotically approach $\gamma_1$ or $\gamma_2$ within $H_2$. \\

Bounded and unbounded motion are defined as follows, \\

\noindent Definition (bounded motion starting on $P_{Y, \Lambda}$) \ \ A point $p \in P_{Y, \Lambda}  \cap H_2$ gives bounded motion for a trajectory $\bm{\psi}(t)$ having $p$ as an initial value at $t=t_0$, if for $t >t_0$, $\bm{\psi}(t)$ does not permanently escape $H_2$ for all $t > t_0$. \\  

It is noted that this definition implies $\bm{\psi}(t)$ could escape $H_2$ at some time $t_1$, but it returns to $H_2$ at a later time $t_2 > t_1$.\\


Both $\hat{S}, W'$ give rise to bounded motion as initial conditions for trajectories in $H_2$.\\


\noindent Definition (unbounded motion starting on $P_{Y, \Lambda}$) \ \ A point $p \in P_{Y, \Lambda} \cap H_2$ leads to unbounded motion for a trajectory $\bm{\psi}(t)$ having $p$ as an initial value at $t=t_0$, if for $t >t_0$, $\bm{\psi}(t)$  moves into $H_1$ or $H_O$ and does not return to $H_2$, or if it asymptotically approaches $\gamma_1$ or $\gamma_2$.\\

Unbounded motion could occur for initial points on $\bar{U}_n$,  but points in these sets could also give rise to bounded motion.

\begin{prpty}
Points in $W', \hat{S}$ give bounded motion. Points in  $\bar{U}_n$ may lead to unbounded motion.
\label{prpty:8}
\end{prpty}

The proof follows by definition of bounded and unbounded motion, together with the dynamics of trajectories with initial values on   $W', \hat{S}$  and $\bar{U}_n$.\\

Define $M^* = W' \cup \hat{S}$, \ $\bar{M}^* = \bar{U}_n$.   This implies,

\begin{prpty}
If $p \in  M^*$  then $p$ yields bounded motion as an initial value for a trajectory.  If $p \in \bar{M}^*$, then $p$ may lead to unbounded motion.
\label{prpty:9}
\end{prpty}

\begin{prpty}
$P_{Y, \Lambda} = M^* \cup \bar{M}^*$. \   $M^* \cap \bar{M}^* = \emptyset $.
\label{prpty:10}
\end{prpty}

The definitions of bounded and unbounded motion for trajectories starting on a point in $P_{Y, \Lambda} \cap H_2$ can be similarly made for iterations of $\Phi(p)$ on $S_{\theta}$ where $p \in S_{\theta} \cap \Lambda$.  Let $\Phi_r(p)$ be the projection of $\Phi(p)$ onto the $r$-coordinate. 
\medskip

Let $\rho$ = distance from $P_2$ to $L_2$. Assume $C \in I$. 
\medskip

$\Phi^n(p)$ is unbounded for $n \geq 1$  if $|\Phi_r^n(p)| > \rho$ for some $n$, or if there exists an $N$ such that for $n > N$ there are no more iterates on $S_{\theta} \cap H_2$.  If $|\Phi_r^n(p)| < \rho$ for a sequence $n = n_j \geq 1$, $n_{j+1} > n_{j}$, $j = 1,2, \ldots $,  then $\Phi^n(p)$ is bounded.

\begin{sumry}
\noindent
(Properties of $M^*$)  \  \  Let $p \in M^*$, \ $w \in \bar{M}^*$, $C \in I$, then \medskip

\noindent  A.) $\Phi^n (p)$ is bounded for all $n \geq 1$  (for each $\theta$) \\
\noindent  B.) $\Phi^n (w)$ may be unbounded for $n$ sufficiently large, \\
\noindent  C.) $\partial M^* \equiv W'$ is an infinite union of Cantor sets,  over all $\theta \in [0,2\pi], C \in I_C$. It is Cantor set (fractal) for each fixed $\theta, C$.\\
\noindent  D.) $M^*$ is compact,\\
\noindent  E.) $P_{Y, \Lambda} = M^* \cup \bar{M}^*$, \   $M^* \cap \bar{M}^* = \emptyset $.\\
\noindent  F.) $M^*$ is defined in the initial parameter plane $P_{Y, \Lambda}$ for $\Phi$ and $\Phi$ is iterated on the iteration plane $S_{\theta}$.
\label{sumry:1}
\end{sumry}

The proof follows by previous results.  \\

In summary, 
\medskip

\begin{thm}
The boundary of $M^*$,  is $W'$, an infinite union of Cantor sets, over all $\theta \in [0, 2\pi], C \in I_C$. It defines  the weak stability boundary, for unstable initial points that give  infinite cycling of trajectories about $P_2$.   The interior of $M^*$ is the set $\hat{S}$ of initial points for trajectories that cycle about $P_2$ infinitely many times in a stable manner, over all $\theta \in [0, 2\pi], C\in \hat{I}_C$. \ $\mu \in I_{\mu}$ is fixed. 
\label{thm:M*W'}
\end{thm}
\medskip\medskip

It is interesting to remark that the properties of $M^*$ are similar to those of the classical Mandelbrot set. This is only given since it seems interesting to note. The definition and properties of a Mandelbrot set, $M$,  
are summarized \cite{Mandelbrot:1980, Devaney:2022}.   The following description of this remark can be skipped since it is not used in the analysis. \\ 

$M$ is defined for the complex map, $Q_c (z) = z^2 + c$, $z \in \mathbb{C}$, and $c \in \mathbb{C}$ a parameter.  $M$ is defined in the $c$-plane and iterates in $z$-plane. 
The iterates are considered for the critical point $z=0$.  That is,  $Q_c^n(0)$, $n \geq 1$.

More precisely, $M$ is defined by those values of $c$ such that the iterates $Q_c^n(0)$ in the $z$-plane are bounded for $n \geq 1$. 
$M$ has the properties,

\begin{sumry}
(Properties of $M$) \ \   Let $c \in M$, \ $\hat{c} \in \bar{M}$, 
\medskip

\noindent  A.) $Q_c^n(0)$ are bounded for all $n \geq 1$,\\
\noindent  B.)  $Q_{\hat{c}}^n(0)$ is unbounded for $n \rightarrow \infty$, \\
\noindent  C.) $\partial M$  is a connected fractal set, \\
\noindent  D.) $M$ is compact,\\
\noindent  E.) $\mathbb{C}  = M\cup \bar{M}$, \   $M\cap \bar{M} = \emptyset $,\\
\noindent  F.) $M$ is defined in the parameter plane $c \in \mathbb{C}$ for $Q_c (0)$  and $Q_c (0)$ is iterated on the complex $z$-plane.
\label{sumry:2}
\end{sumry}

\begin{sumry}
 (Differences between $M^*, \  M$)   \ \  
\medskip

\noindent  i.) The iterates $Q^n_c(0)$ on the $z$-plane are parameterized for different $c \in \mathbb{C}$;  The iterates $\Phi^n(c)$ on the $(r,\dot{r})$-plane are parameterized different $c \in P_{Y, \Lambda} \subset \mathbb{R}^2$\medskip

\noindent  ii.)   $\partial M $ is a connected fractal curve; $ \partial M^* = W'$ is fractal for each fixed $\theta, C$, and totally disconnected. \medskip

\noindent  iii.)   For points in $\bar{M}$, the iterates are unbounded, and for the points in $\bar{M}^*$, the iterates may be unbounded.\\ 
 
\label{sumry:3}
\end{sumry}

$M^*$ and $M$ have some similarities. They are both defined in a parameter plane ($M^*$ in  $P_{Y, \Lambda}$ and $M$ in the $c$-plane), and their iterates in another plane. They both have fractal boundaries, with the main difference that $\partial M^*$ is totally disconnected for each $\theta, C$, whereas $\partial M$ is connected. $M$ is defined for a complex analytic map, whereas, $M^*$ is defined for a real analytic map. This is summarized in Table \ref{tab:Comparison}.  \\

\begin{table}[]
\begin{tabular}{|c|c| }
\hline
$M^*$ & $M$\\ 
\hline   \\
$\Phi^n(c), n\geq 1, $ bounded, $c \in P_{Y, \Lambda} \subset \mathbb{R}^2$ &  $Q_c^n(0)$, $n \geq 1$,  bounded $c \in \mathbb{C}$  \\
& \\
$\partial M^*$ fractal, totally disconnected for each $\theta,C$ & $\partial M$ fractal, connected\\
& \\
$\Phi^n(c) \in S_{\theta} \subset \mathbb{R}^2$ &    $Q_c^n(0)  \in  \mathbb{C}$\\
\hline 
\end{tabular}
\medskip\medskip

\caption{Comparison between $M^*$ and $M$.}
\label{tab:Comparison}
\end{table}

\section{Non-existence of KAM tori}
\label{sec:NoKAM}
\medskip

KAM tori exist for $C$ sufficiently large as follows from \cite{Kummer:1979}.  See Theorem \ref{thm:KAM} for $C \in C_L$. When $C \notin C_L$  then the general existence of KAM tori in $H_2$ about $P_2$ is an open problem to theoretically understand.  When $C$ is not large and the radius of $H_2$ increases, then the gravitational perturbative force from $P_1$ can be strong enough relative to the gravitational force of $P_2$, for $\mu$ sufficiently small, so that KAM tori may not exist on $\Sigma(C)$.

The range of $C$ considered is $C \in I$,  and $\mu \in I_{\mu}$ is fixed. The techniques in  \cite{Kummer:1979} are  not applicable for $C \in I$.

  \begin{thm}
Consider a  fixed energy surface $\Sigma(C), \ C \in I_C$, $\mu \in I_{\mu}$,  and a section $S_{\theta}$ on this energy surface for any given $\theta \in [0, 2\pi]$. The iterates of an 
initial point $p_0 \in P_{Y, \Lambda}$, $\Phi^k(p_0)$ on $S_{\theta}$, cannot lie on a KAM torus if $p_0 \in \partial M^*$.  If $p_0 \in \bar{M}^*$ the iterates may not lie on a KAM torus.   If $p_0 \in \text{the interior of} \ M^*$, then it is possible the iterates could lie on a KAM torus.
\label{thm:NoKAM}
\end{thm}

 \noindent Proof  -

Consider a trajectory ${\bm{\psi}}(t,\mu)$ cycling about $P_2$ in ${\bf Y}$-coordinates under these assumptions for $t > 0$, with an inital value at $t=t_0$ on $M^* \subset P_{Y, \Lambda}$.


The points of intersection of ${\bm{\psi}}(t,\mu)$ with $S_{\theta} \cap \Sigma(C)$ occur at a sequence of times $t_k, k\geq 1 $, $t_{k+1} > t_k$. These generate the iterates of $\Phi$, labeled $\Phi^{k}(p_0)$, where $p_0$ is the initial point of ${\bm{\psi}}(t,\mu)$ on $S_{\theta}$ at $t=t_0$.

$\Phi^{k}(p_0)$, $k=1,2, \ldots$ lie on a Cantor set $\mathcal{C}_C$  on $S_{\theta}$ by assumptions.

 The $\Phi^{k}(p_0)$,\ $k=1,2, \ldots$  cannot lie on a KAM torus.  This is shown as follows:   Assume the iterates did lie on such a KAM torus, $T^*$, two-dimensional. It  intersects $S_{\theta}$ in a topological circle, $\tilde{S}^1$.  Since the iterates belong to  $\mathcal{C}_C$, they are nowhere dense.  This contradicts the assumption they lie on $\tilde{S}^1$, since then the iterates would have to be dense on $\tilde{S}^1$ by the Moser Twist Theorem for a monotone twist map, $M$ \cite{Moser:1962, SiegelMoser:1971}.  This follows since by KAM theory,  $M$ takes the form,  $\phi \rightarrow \phi + 2\pi \lambda(\rho) + \mathcal{O}(\mu), \ \rho \rightarrow \rho + \mathcal{O}(\mu)$, $\phi \in [0, 2\pi], \rho > 0$ are polar coordinates for points on $\tilde{S}^1$, where the angular frequency for a given value of $\rho$, $\lambda$, satisfies, $d \lambda/d\rho > 0$, and satisfies diophantine conditions. Thus the iterates cannot lie on $T^*$.
   
 These iterates, however, could belong to points in a resonance gap between KAM tori, if they existed, on $\Sigma(C)$.   
If that were the case, then how wide could this gap be on $\Sigma(C)$? In the complementory region $\bar{M}^* \cap \Sigma(C)$,  initial values for $\Phi$ may not yield iterates that lie on the intersection of KAM tori with $S_{\theta} \cap \Sigma(C)$. This is because the region $\bar{M}^*$ consists of points in $\bar{U}_n$ for finite $n \geq 1$. They may be unbounded by iteration with $\Phi$ on $S_{\theta}$.   Thus, this indicates that initial points in $\bar{M}^*$ to use for iteration of $\Phi$ on $S_{\theta} \cap \Sigma(C)$ may not lie on KAM tori. 

What happens for $p_0 \in \text{the interior of} \ M^*$? These points belong to $\hat{S}_{C}(\theta)$, cycling about $P_2$ for all time,  and could lie on KAM tori.  This is not known.

This proves Theroem \ref{thm:NoKAM} \\ \\

This theorem only pertains to the fate of the iterations of points on $P_{Y, \Lambda}$.

\section{Discussion of Results}
 \label{sec:Discussion}
\medskip

One of the main results of this paper is that the weak stability boundary, $W'$, about $P_2$ for infinitely many cycles of $P$ about $P_2$, 
is a union of infinitely many Cantor sets of hyperbolic points.     This answers the question on the fractal structure of the weak stability boundary, at least in the 
case of infinite cycling. { However, this result is not competely rigorous and relies on some semi-analytic and numerical results used in several assumptions.}

Another result is that KAM tori cannot exist on $W'$, but they may exist on the stable set $\hat{S}$ that it bounds for the range of $C$. Results indicate that that KAM tori may not exist beyond $W'$, so this set may yield a kind of boundary for KAM tori about $P_2$. 

A curious result is that $M^* = W' \cup \hat{S}$  has similar properties to a classical Mandelbrot set, even though they are defined completely differently. The main difference between $M^*$ and $M$ is that although the boundaries are both fractal sets, the boundary of $M$ is continuous whereas the boundary of $M^*$ is totally discontinuous.

This paper uses previous numerical results to motivate the analysis.  This is done for relatively small values of $n$.  It would be interesting to numerically explore $M^*$ for large values of $n$, which is beyond the scope of this paper.

 A potential application of these results pertains to the low energy permanent weak capture of $P$ about $P_2$ from the Hill's region about $P_1$. $P$ would move through the channel, $\mathcal{C}$, and into the Hill's region about $P_2$ to the boundary region, $W'$ about $P_2$ and become permanetly captured taking infinitely many cycles.
The location of  $W'$ relative to $P_2$ would give regions where permanent capture can occur. Permanent capture using $W'$ could be used to design spacecraft trajectories 
that never need orbit maintenance maneuvers while orbiting $P_2$. 
   
The results of this paper may have an interesting interpretation on the nature of  the interaction of the gravitational fields between two bodies on a particle of negligible mass.
Defining this boundary over infinite cycles yields a boundary consisting of the union of Cantor sets.   Thus, in general,  no matter how small one magnifies this boundary, the self similarity yields the same structure. This could just be a curious mathematical property or it also may say something more about interacting gravitational fields and their scale properties.

\appendix
\section{Supporting Calculations}
\label{sec:Appendix} 

\subsection{Kepler Energy,  $\bm{E_2}$}
\label{subsec:KeplerEnergy}

In a $P_2$-centered inertial coordinate system,  ${\bf{X} }= (X_1, X_2)$,  the Kepler energy of $P$ relative to $P_2$ is
\begin{equation}
E_2 =  (1/2)|\dot{\bf{X}}|^2 - \mu |{\bf X}|^{-1},
\label{eq:KepEnergyInertial}
\end{equation}
where $X_1 = x_1 -1, X_2 = x_2$.  ${\bf{x}} = (x_1, x_2)$ are $P_1$-centered inertial coordinates. .
\medskip

In a $P_2$-centered rotating coordinates,  ${\bf{Y}} = (Y_1, Y_2)$, obtained by setting $Y_1 = y_1 - 1, Y_2= y_2$, where ${\bf{y}} = (y_1, y_2)$ are $P_1$-centered rotating coordinates defined
in (\ref{eq:DEs});
\begin{equation}
E_2({\bf Y}, \dot{\bf{Y}}) = (1/2) \dot{Y}^2 - \frac{\mu}{Y} -L({\bf Y}, \dot{\bf{Y}}) + (1/2) Y^2,
\label{eq:KepEnergyRot}
\end{equation}
where $L({\bf Y}, {\bm{\dot{Y}}}) = \dot{Y}_1Y_2 - \dot{Y}_2 Y_1$, $Y = |{\bf{Y}}|, \dot{Y} = |\dot{\bf{Y}}|$.

\subsection{$\bm{E_2 \geq 0}$ When $\bm{P}$ Cycles About $\bm{P_1}$ }
\label{subsec:UnstableCyclingMotionPosEng}

It is assumed $C_2 < C\lessapprox C_1$.  $P$ has an initial conditon at $t=t_0$ on $L(\theta)$ with eccentricity $e \in [0,1)$, and with initial velocity (\ref{eq:vp}),  at a distance $r$ from $S$. This implies at the initial point,  $E_2 < 0$. $P$ cycles around $S$ $n \geq 1$ times and escapes $S$ through the channel $\mathcal{C}_1$.  Since this an unstable motion, then this implies the trajectory for $P$ will make at least one cycle about $P_1$. In particular, $P$ will cross the  negative $Y_1$-axis, to the left of $P_1$, where $r > 1$. It will do this while moving near the zero velocity curve about $P_1$ (see (Belbruno2004, Figure 3.10), (McGehee1969)).  Thus, there exists a time $t_ 1 > t_0$, where $Y_2(t_1) =0, Y_1(t_1) \gtrapprox -2$.  Also, for the velocities, $|\dot{Y}_2(t_1)|> 0$ since $P$ will cross the $Y_1$-axis moving upward or downward, and $\dot{Y}_1(t_1) \approx 0$ since $P$ crosses the 
$Y_2$-axis near the zero velocity curve.  Thus, at the time of crossing, $L(t_1)  \approx 0$.   This implies $E_2(t_1) \approx (1/2) \dot{Y}^2 +  (1/2) Y^2 - \frac{\mu}{Y}$ which is positive since $\mu$ is small. 
Thus, $E_2(t_1) > 0$. 

\subsection{Definition of  $\bm{S_{\theta_0}}$  }
\label{subsec:SolvingforThetaDot}

In a $P_2$-centered rotating system,  $(Y_1, Y_2)$, the Jacobi integral in the $S$-centered coordinates is 
\begin{equation}
\tilde{J} = -(|\dot{\bf{Y}}|^2)  + 2 [ (1-\mu)/r_1   +   \mu/r ]+  [(Y_1 + 1 -\mu)^2 + Y_2^2]  + \mu(1-\mu),
\label{eq:JacobiSCentered}
\end{equation}
where   $r_1^2 = (Y_1-1)^2 + Y_2^2$, $r = |{\bf{Y}}|$.  $\tilde{J} = C$ defines the energy surface, $\Sigma  \Sigma(C)$.

In polar coordinates,  $(r, \theta)$, where $Y_1 = r\cos\theta, Y_2 = r\sin \theta$, it is verified that the Jacobi integral, $\tilde{J}$ becomes
\begin{equation}
\hat{J} =   -(\dot{r}^2 + r^2 \dot{\theta}^2) + f(r, \theta),
\label{eq:JacobiPolar}
\end{equation}
where $f$ is well defined for $r >0, r_1 > 0$, and does not depend on $\dot{r}, \dot{\theta}$.  Thus, $\hat{J} =C$, implies that 
$\dot{\theta}^2 = \dot{\theta}^2 (r, \dot{r}, \theta) = r^{-2}[\dot{r}^2 + f(r, \theta)]$.  Choosing the positive root, 
$\dot{\theta} = \dot{\theta} > 0$, as required for $S_{\theta_0}$ on $\hat{J}=C$, $S_{\theta_0} = \{r, \dot{r}| \theta=\theta_0,\}$, where  $\dot{\theta}= \dot{\theta}(r, \dot{r}, \theta)$. \\

\subsection{Trajectory Properties} 
\label{subsec:TrajProp}
 
The trajectories satisfy topological conditions.  The range of $C$ is $C_a < C < C_2$ and $C_2 \leq C < C_1$, These conditions are:
 \medskip

\noindent {\it Hypothesis $A$}
\medskip

\noindent
(i) All the trajectories on $W^s(\gamma_1)$ in $H_2$ make at least $n$ consecutive cycles about $P_2$, $n=1,2, \ldots$, \\

\noindent (ii) All the trajectories  on $W^u(\gamma_1)$ in $H_1$ must make at least $1$ cycle about $P_1$,\\

\noindent (iii) All the trajectories on $W^s(\gamma_2)$ in $H_2$  must make at least $n$ consecutive cycles about $P_2$, \\

\noindent (iv) All the trajectories  on $W^u(\gamma_2)$ in $H_O$ must make at least $1$ cycle about $P_1$.\\

In Hypothesis $A$, cases occur where these assumptions are not satisfied. 
For this and other details, see\cite{BelbrunoGideaTopputo:2010}. {This hypothesis is discussed in Assumption \ref{assum:3}.}\\

\noindent
{\bf Acknowledgements}
\medskip

\noindent
 Many thanks to Marian Gidea for helpful comments. I would like to thank Princeton University, Dept. of Astrophysical Sciences.   Research partially funded by NSF grant DMS-1814543.  
\medskip



\begin{thebibliography}{99}

\bibitem{Arnold:1989}
 Arnold, V. I.:
 Mathematical Methods of Classical Mechanics.
 GTM {\bf 60} , Springer Verlag, Berlin   (1989)


\bibitem{Belbruno:1987}
 Belbruno,  E.:
 Lunar capture orbits,  a method of constructing Earth-Moon trajectories and the Lunar GAS mission, with Ballistic Capture and the Lunar GAS Mission.
 in Proceedings of the 19th AIAA/DGLR/JSASS Inter. Elec. Propl. Conf.,  Paper no. 87-1054, Colorado Springs   (1987)  (arc.aiaa.org/doi/10.2514/6.1987-1054).
  

\bibitem{BelbrunoMiller:1993}
Belbruno, E., Miller, J.: 
Sun-perturbed Earth-to-Moon transfers with ballistic capture.
J. Guidance, Control, Dynamics {\bf 16}, 770-775  (1993)

\bibitem{BelbrunoMarsden:1997}
Belbruno, E.,  Marsden, B.: 
Resonance hopping in comets.
AJ  {\bf  113},  1433-1444  (1997)


\bibitem{Belbruno:2004}
 Belbruno, E.: 
Capture Dynamics and Chaotic Motions in Celestial Mechanics.
 Princeton University Press, Princeton   (2004)
    
\bibitem{Belbruno:2007}
Belbruno, E.:
Fly Me to the Moon: An Insiders Guide to the New Science of Space Travel.
Princeton University Press, Princeton   (2007)

\bibitem{BelbrunoGideaTopputo:2010}
Belbruno, E.,  Gidea, M., Topputo, F.   
Weak stability boundary and invariant manifolds.
SIAM J Appl. Dyn. Sys., {\bf  9},   1061-1089   (2010)
  


\bibitem{Celletti:2007}
Celletti, A.,  Chierchia, L. 
KAM stability and celestial mechanics.
Memoirs of the American Mathematical Society {\bf 187}   (2007)

\bibitem{Circi:2001}
Circi, C., P. Teofilato, P.:
On the dynamics of weak stability boundary lunar transfers.
Cel. Mech. Dyn. Astr. {\bf 79}, 41-72  (2001)


\bibitem{Conley:1968}
Conley, C.: 
Low energy transit orbits in the restricted three-body problem.
SIAM J. Appl. Math.  {\bf 16},   732-746  (1968)

\bibitem{Conley:1969}
Conley, C.:
On the ultimate behavior of orbits with respect to an unstable critical point.
J. Diff Equ, {\bf  5}, 136-158  (1969)

\bibitem{Devaney:2022}
Devaney, R.:
A First Course in Chaotic Dynamical Systems.
  CRC Press, Boca Raton, Florida (2022)




\bibitem{GarciaGomez:2007}
 Garcia, F., Gomez, G.: 
A Note on Weak Stability Boundaries.
Cel. Mech. Dyn. Astr. {\bf 97}, 87-100 (2007)
  
\bibitem{GideaMasdemont:2007}
Gidea, M.,  Masdemont, J.: 
Geometry of homoclinic connections in the planar circular restricted three-body problem.
 Inter. J. Bifur. Chaos. Sci. Engrg. {\bf 17}, 1151-1169   (2007)


\bibitem{GuckenheimerHolmes:1983}
Guckenheimer, J.,  Holmes, P.:
Nonlinear Oscillations, Dynamical Systems, and Bifurcations of Vector Fields.
  Applied Mathematical Sciences {\bf 42}
Springer Verlag, Berlin (1983) 

\bibitem{Jehn:2008}
Jehn, R.,  et. al.:
Navigating BepiColombo during the weak-stability capture at Mercury.
Advances in Space Research  {\bf 42}, 1364-1369 (2008)

\bibitem{JianYiSui:2015}
 Jian, L.,  Yi-Sui, S.: 
A survey of weak stability boundaries in the Sun-Mars system.
RAA {\bf 15} 376 (2015)



\bibitem{KoonLoMarsdenRoss:2000}
Koon, W.S., Lo, M.W.,  Marsden, J.E.,  Ross, S.D.: 
Heteroclinic connections between periodic orbits and resonance transitions in celestial mechanics.
Chaos  {\bf 10} ,  427-469  (2000)


\bibitem{Kummer:1979}
Kummer, M.:
 On the stability of Hill's solutions of the plane restricted three-body problem.
Am. J. Math. {\bf  101}, 1333-135 (1979)

\bibitem{LlibreMartinezSimo:1985}
Llibre, J., Martinez, R., Sim\'o, C.: 
Transversality of the invariant manifolds associated to the Lyapunov family of periodic orbits near $L_2$ in the
Restricted three-body problem.
J. Diff Equ {\bf 58}  104-156 (1985)





\bibitem{Mandelbrot:1980}
Mandelbrot, B.:
Fractal aspects of the iteration of $z \mapsto \lambda z(1-z)$ for complex $\lambda, z^n$,
Annals of the New York Academy of Sciences {\bf 357},  249-259 (1980)


\bibitem{Milnor:2006}
Milnor, J.:
Dynamics in One Complex Variable. 
 Annals of Mathematics Studies,  160, third edition,
Princeton University Press, Princeton  (2006)


\bibitem{Moser:1962}
Moser, J.:
On invariant curves of area-preserving mappings of an annulus.
Nachrichten der Akademie der Wissenschaften Gottingen, 
Math.-Phys. Kl. {\bf 2}, 1-20 (1962)


\bibitem{Moser:1973}
 Moser, J.:
 Stable and Random Motions in Dynamical Systems.
  Princeton University Press, Princeton  (1973)




\bibitem{Poincare:1899}
 Poincar\'{e}, H.:
 Les Methodes Nouvelles de la M\'{e}chanique Celeste, Vols 1-3.
  Gauthier Villars, Paris  (1899)

\bibitem{Romagnoli:2009}
Romagnoli, D., Circi, C.:
Earth-Moon weak stability boundaries in the restricted three and four body problem.
Cel. Mech. Dyn. Astr. {\bf 10}, 79-103 (2009)



\bibitem{Schoenmaekers:2001}
Schoenmaekers, J., Horas, D., Pulido, D.: 
 SMART-1: With solar electric propulsion to the Moon.
in Proceedings of the 16th International Symposium on Space Flight Dynamics, Pasadena, CA, 114-130 (2001)

\bibitem{SiegelMoser:1971}
 Siegel, C.L., Moser, J.K.:
  Lectures on Celestial Mechanics.
 Grundlehren Series {\bf 187}, Springer-Verlag, Berlin  (1971)


\bibitem{Simo:2000}
Sim\'o, C., Stuchi, T.: 
Central stable/unstable manifolds and the destruction of KAM tori in the planar Hill problem.
Physica D {\bf 140},  1-32.  (2000)



\bibitem{TopputoBelbruno:2009}
Topputo, F., Belbruno, E.:  (2009),
Computation of weak stability boundaries: Sun-Jupiter. 
Cel. Mech. Dyn. Astr., {\bf 105}, 3-17  (2009)



\end{thebibliography}
\end{document}